\documentclass{scrartcl}

\usepackage{lineno,hyperref}
\modulolinenumbers[5]

\usepackage[utf8]{inputenc}
\usepackage[english]{babel}
\usepackage{amsmath,amssymb,upref}
\usepackage{amsfonts}
\usepackage{amsthm,array,makecell}
\usepackage{color}
\usepackage{graphicx}
\usepackage{sidecap}
\usepackage{multirow}
\usepackage{booktabs}
\usepackage{fullpage}

\usepackage[square,sort,comma,numbers]{natbib}

\usepackage{caption}
\usepackage{subcaption}

\usepackage{longtable}
\usepackage{rotating}

\usepackage{algorithm}
\usepackage{algpseudocode}

\usepackage{pgfpages}
\usepackage{tikz}
\usepackage{grffile}
\graphicspath{{figures/}}

\usetikzlibrary{decorations.text}
\usetikzlibrary{shapes,snakes}
\usetikzlibrary{shapes}
\usetikzlibrary{arrows}
\usetikzlibrary{positioning}
\usetikzlibrary{matrix}
\usetikzlibrary{patterns}
\usetikzlibrary{decorations.pathreplacing,decorations.pathmorphing,decorations.markings}
\usetikzlibrary{positioning,tikzmark}
\usepackage{pgfplots}

\usepackage{xcolor}
\definecolor{dark-violet}{RGB}{148,0,211}
\definecolor{sea-green}{RGB}{46, 139,  87}

\DeclareGraphicsExtensions{.pdf,.eps,.png,.jpg,.jpeg}

\usepackage{listings}
\usepackage{mathrsfs}
\usepackage{color}
\usepackage{amssymb}
\usepackage[utf8]{inputenc}
\usepackage{afterpage}
\usepackage[T1]{fontenc}
\usepackage{sectsty}
\usepackage{amsfonts}
\usepackage{mathtools}
\mathtoolsset{showonlyrefs}
\usepackage{amsmath}
\usepackage{mathrsfs}
\usepackage{amsthm}
\usepackage{verbatim}
\usepackage[margin=1in]{geometry}
\usepackage{bm}
\usepackage{bbm}
\usepackage{hyperref}
\usepackage{enumitem}
\usepackage{natbib}
\usepackage{lastpage}
\usepackage{hyperref}

\graphicspath{{figures/}}

\newcommand*\mcupinn[2]{\vcenter{\hbox{$\mathsurround=0pt
  \ifx\displaystyle#1\textstyle\else#1\fi\bigcup$}}}

\newcommand*\mcapinn[2]{\vcenter{\hbox{$\mathsurround=0pt
  \ifx\displaystyle#1\textstyle\else#1\fi\bigcap$}}}

\newcommand{\R}{\mathbb{R}}

\newcommand{\hu}{\hat{u}}

\newcommand {\bx} {{\boldsymbol{x}}}

\newcommand {\ba} {{\boldsymbol{a}}}

\newcommand {\bc} {{\boldsymbol{c}}}

\newcommand{\Xcal}{\mathcal{X}}

\newcommand{\Ucal}{\mathcal{U}}
\newcommand{\nh}{N}
\newcommand{\hbfM}{\hat{\bfM}}
\newcommand{\hbfF}{\hat{\bfF}}

\newcounter{matriz}

\newcommand{\bftheta}{\boldsymbol{\theta}}
\newcommand{\bfx}{\boldsymbol{x}}

\newcommand{\bfmu}{\boldsymbol{\mu}}
\newcommand{\bfSigma}{\boldsymbol{\Sigma}}

\newcommand{\bfM}{\boldsymbol{M}}
\newcommand{\bfF}{\boldsymbol{F}}

\newcommand{\nm}{m}

\newcommand{\dbftheta}{\dot{\bftheta}}
\newcommand{\bfy}{\boldsymbol y}

\ifpdf
\hypersetup{
  pdftitle={Coupling parameter and particle dynamics for adaptive sampling in Neural Galerkin schemes}, 
  pdfauthor={Yuxiao Wen and Eric Vanden-Eijnden and Benjamin Peherstorfer} 
}
\fi

   {\begin{trivlist}\item[]{\bfseries\sffamily Keywords:}\ }
   {\end{trivlist}}
   
   \title{Coupling parameter and particle dynamics for adaptive sampling in Neural Galerkin schemes}

\author{Yuxiao Wen\footnote{Courant Institute of Mathematical Sciences, New York University, 251 Mercer Street, New York, NY 10012, USA} \and Eric Vanden-Eijnden\footnotemark[1] \and Benjamin Peherstorfer\footnotemark[1]}

\begin{document}

\maketitle

\begin{abstract}
Training nonlinear parametrizations such as deep neural networks to numerically approximate solutions of partial differential equations is often based on minimizing a loss that includes the residual, which is analytically available in limited settings only. At the same time, empirically estimating the training loss is challenging because residuals and related quantities can have high variance, especially for transport-dominated and high-dimensional problems that exhibit local features such as waves and coherent structures. Thus, estimators based on data samples from un-informed, uniform distributions are inefficient. This work introduces Neural Galerkin schemes that estimate the training loss with data from adaptive distributions, which are empirically represented via ensembles of particles. The ensembles are actively adapted by evolving the particles with dynamics coupled to the nonlinear parametrizations of the solution fields so that the ensembles remain informative for estimating the training loss. Numerical experiments indicate that few dynamic particles are sufficient for obtaining accurate empirical estimates of the training loss, even for problems with local features and with high-dimensional spatial domains.
\end{abstract}

\section{Introduction}\label{sec:Intro}
Nonlinear parametrizations (discretizations) that approximate solution fields of partial differential equations (PDEs) on manifolds can achieve faster error decays than traditional, linear parametrizations that are restricted to approximations in vector spaces. A widely used approach for numerically fitting nonlinear parametrizations is minimizing the norm of the PDE residual, oftentimes even globally over the time interval of interest; see \cite{ Weinan2017TheDR,Sirignano_2018, PINN19,Berg_2018,LI2020109338,khoo_lu_ying_2021} and the early works \cite{MN1994,doi:10.1080/00986449208936084,366006}.
However, obtaining an estimate of the norm of the residual---the empirical risk---to find parameters that minimize it can be challenging \cite{pmlr-v145-rotskoff22a,WeinanEOpinion,NG22}. For example, consider a PDE that describes the transport of localized features such as phase transitions, waves, or other coherent structures across the spatial domain. For such transport-dominated problems it has been shown that the faster error decay of nonlinear parametrizations can be beneficial compared to linear approximations from an approximation-theoretic perspective \cite{OHLBERGER2013901,Greif19,10.1093/imanum/dru066,MADAY2002289,Cohen2021,P22AMS}. However, precisely for such problems, an un-informed, passive sampling with a uniform distribution over the spatial domain leads to estimators of the residual norm with high variance because the local features of the solution fields typically lead to localized residuals so that few samples fall in regions of the spatial domain where the magnitude of this residual is large and, thereby, informative. Equations formulated over high-dimensional spatial domains lead to similar challenges of estimating residual norms and training losses in general. The described sampling issue is a special case of the pervasive challenge in machine learning of estimating the population loss with the empirical loss from data \cite{NIPS1991_ff4d5fbb}.

The challenge of accurately estimating the training loss when fitting nonlinear parametrizations has been recognized in a range of works. For certain limited classes of PDEs and corresponding parametrizations, quantities that are needed for computing the training loss can be derived analytically \cite{Lubich2008,doi:10.1137/21M1415972} and special structure can be exploited \cite{anderson2023fast}.  
In model reduction \cite{RozzaPateraSurvey,SIREVSurvey,ManzoniBook2016}, sampling via empirical interpolation \cite{Everson1995,Barrault2004,DEIM,doi:10.1137/15M1019271,PDG18ODEIM,NEGRI2015431} is an important component of reduced models when they are formulated over nonlinear parametrizations that represent latent manifolds rather than subspaces; see \cite{P22AMS} for a survey. Various sampling methods for empirical interpolation and related techniques are proposed in \cite{doi:10.1137/19M1257275,9147832,UngerTransformModes2020,CohenSampling} and the work \cite{MOJGANI2023115810} investigates the sample complexity in the context of nonlinear model reduction. There also is work on using deep autoencoders for dimensionality reduction, where sparse sampling methods are proposed for efficient computations \cite{KIM2022110841,Romor2023}. However, these sampling schemes aim to find an accurate set of samples in a pre-processing step, instead of actively adapting samples.
In the context of high-dimensional problems, many different active learning and active data acquisition approaches to adaptively sample the temporal and spatial domains have been developed. There are works that generate a batch of uniform samples and sub-select and use the ones with highest absolute residual values for training; see, e.g., \cite{Lu_2021,Gao_2021,Wight_2020} and the survey \cite{WU2023115671}. However, the initial batch of samples has to capture already the high-residual regions of the spatial domain, which can be challenging for problems with local features in high-dimensional spatial domains. 
The work \cite{Tang2021} trains a separate normalizing flow for generating samples with high residual values, which in turn requires informative training data. In \cite{GuoWang22}, the time domain is split into subdomains and samples in each subdomain are drawn proportionally to its residual value, which is closely related to domain decomposition techniques.
Other techniques leverage certain structures in the class of PDEs and setups of interest, such as parameterizing the gradient of the PDE solution rather than directly the solution field \cite{E_2017,Han_2018}, approximating committor functions \cite{Khoo2018,doi:10.1063/1.5110439,pmlr-v145-rotskoff22a}, and connections to control problems \cite{E2017,JMLR:v18:17-653,LiControl2023}. Another line of work builds on Gaussian process regression for approximating solution fields and offers strong error analyses, which can be informative for data collection as well \cite{CHEN2021110668,batlle2023error}.

In this work, we estimate the training loss with an ensemble of particles that we actively adapt by imposing dynamics on the particles  so that they stay informative for estimating the training loss even as the PDE solution field evolves in time. The adaptation of the particles leads to an active data acquisition in the sense that the training loss is measured in a dynamic way so that few measurements and thus samples are sufficient for an accurate estimation. 
The starting point for us is the Dirac-Frenkel variational principle \cite{dirac_1930,Frenkel1934,Lubich2008,lasser_lubich_2020} that we build on to formulate a Galerkin problem \cite{NG22} that is solved sequentially in time; see also \cite{doi:10.1137/050639703,SAPSIS20092347,hesthaven_pagliantini_rozza_2022,ZPW17SIMAXManifold,Du_2021,doi:10.1137/21M1415972,finzi2023a,chen2023crom}.  
Our sequential-in-time approach leads to a dynamical system for the parameters and this dynamical system is integrated forward in time with standard time-integration schemes to compute the parameters that represent the approximate solution field over time.
Fitting the parameters sequentially in time with a dynamical system is key for deriving a dynamical system for the ensemble of particles so that it can be moved along in time and lead to an efficient training-loss estimation. It is important to note that with our approach the particle dynamics can be derived for a wide range of nonlinear parametrizations, which is in stark contrast to the previous work \cite{NG22} that also proposes an active learning scheme but one that is only applicable with a limited class of nonlinear parametrizations, namely shallow neural networks with exponential units. 
Numerical experiments indicate that integrating forward in time the coupled parameter and particle dynamics circumvent the sampling issue even for problems with local features and with high-dimensional spatial domains.

The paper is organized as follows. In Section~\ref{sec:Prelim} we briefly review nonlinear approximations methods based on the Dirac-Frenkel variational principle such as Neural Galerkin schemes and mathematically formulate the challenge of sampling to compute residual-based quantities for training. In Section~\ref{sec:adaptive_samp}, we derive Neural Galerkin schemes with coupled parameter and particle dynamics. Numerical results are reported in Section~\ref{sec:numerical_results} and conclusions are drawn in Section~\ref{sec:Conc}.

\section{Preliminaries}\label{sec:Prelim}
We briefly describe the need for nonlinear parametrizations for approximating efficiently solutions of certain classes of PDEs in Section~\ref{sec:Prelim:Need} and discuss Neural Galerkin schemes based on the Dirac-Frenkel variational principle in Section~\ref{sec:Prelim:NG} for numerically solving PDEs with nonlinear parametrizations. The problem formulation is given in Section~\ref{sec:Prelim:Prob}, where we describe that training nonlinear parametrizations is challenging because solutions of precisely those classes of PDEs where nonlinear parametrizations are beneficial in terms of approximation theory often exhibit features that are local in the spatial domain. 

\subsection{The need for nonlinear parametrizations}\label{sec:Prelim:Need}
Consider the following PDE
\begin{align}
\label{eq:PDE}
\partial_tu(t,\bx) &= f(\bx,u),& (t,\bx)\in[0,\infty)\times\mathcal{X},\\
u(0,\bx) &= u_0(\bx),&\bx\in\mathcal{X},
\end{align}
where $u:[0,\infty)\times\mathcal{X}\rightarrow\mathbb{R}$ is a time-dependent field in a space $\Ucal$ over the spatial domain $\mathcal{X}\subseteq\R^d$  and $u_0:\Xcal\rightarrow\R$ is the initial condition.  The right-hand side function $f: \Xcal \times \Ucal \to \mathbb{R}$ can include partial derivatives of $u$ in the spatial coordinate $\bfx$. Note that $f$ could depend on time $t$ and it is only for ease of exposition that we skip the time dependence in the following. We consider either Dirichlet or periodic boundary conditions to make the problem \eqref{eq:PDE} well posed, which means that the boundary conditions can be imposed with a suitable choice of the space $\Ucal$ so that all functions in $\Ucal$ satisfy the boundary conditions. As an alternative and for other boundary conditions, they can be imposed via penalty terms, as we demonstrate with numerical experiments in Section~\ref{sec:numerical_results}.

A typical approach in scientific computing to numerical solve PDEs such as \eqref{eq:PDE} is to compute an approximate solution $\tilde{u}$ in a finite-dimensional subspace $\Ucal_{\nh}$ of $\Ucal$. The subspace $\Ucal_{\nh}$ is spanned by $\nh$ basis functions $\phi_1, \dots, \phi_{\nh}$ and the approximation $\tilde{u}$ is a linear combination of the basis functions with $\nh$ coefficients  $\theta_1(t), \dots, \theta_{\nh}(t) \in \mathbb{R}$. Such an approximation is linear in the parameter $\bftheta(t) = [\theta_1(t), \dots, \theta_{\nh}(t)]^T \in \mathbb{R}^{\nh}$. There are classes of PDEs where it is known that such linear approximations in subspaces can require high dimensions $\nh$. Often such problems are transport-dominated in the sense that a coherent structure such as a wave is traveling through the spatial domain. Notice that adaptive mesh refinement was introduced for efficiently approximating such problems \cite{BERGER198964,BergerLeVeque98}. We refer to the notion of Kolmogorov $\nh$-width \cite{10.1093/imanum/dru066,MADAY2002289,Cohen2021,OHLBERGER2013901,Greif19,ManzoniNwidth,Nonino2023} and nonlinear model reduction \cite{P22AMS} for more details about such problems. 
Another class of PDEs for which approximations in subspaces as described above becomes challenging is given by PDEs with high-dimensional spatial domains $\Xcal$. Discretizing even moderately high-dimensional spatial domains $\Xcal$ with regular grid-based techniques based on linear approximations as described above becomes intractable quickly because of the curse of dimensionality.

\subsection{Nonlinear parametrizations and the Dirac-Frenkel variational principle} \label{sec:Prelim:NG}
In cases where linear approximations in subspaces lead to slow error decay, one can resort to parametrizations $\hu: \Xcal \times \Theta \to \mathbb{R}$ that are nonlinear in the sense that the parameter $\bftheta(t) \in \Theta \subseteq \mathbb{R}^{\nh}$ enters nonlinearly in the second argument of $\hu$. One example of such a nonlinear parametrization is given by deep neural networks that have time-dependent weights $\bftheta(t)$. In the following, we restrict ourselves to parametrizations $\hu$ that also impose the boundary conditions that accompany the PDE in \eqref{eq:PDE}.

We now follow \cite{NG22} and build on the Dirac-Frenkel variational principle \cite{dirac_1930,Frenkel1934,Lubich2008} to derive a dynamical system for $\bftheta(t)$ such that the corresponding parametrization $\hu(\cdot; \bftheta(t))$ numerically solves the PDE in \eqref{eq:PDE} over time; see also \cite{doi:10.1137/050639703,SAPSIS20092347,Du_2021,doi:10.1137/21M1415972,lasser_lubich_2020} for related work. The process is analogous to the classical method of lines approach for numerically solving PDEs: The PDE \eqref{eq:PDE} is first discretized in space, which results in a dynamical system---system of ordinary differential equations---that is then discretized and integrated forward in time to compute a numerical solution of the PDE. 
Let us start with the residual of the PDE \eqref{eq:PDE} at time $t$
\begin{equation}
r_t(\bx; \bftheta(t),\dbftheta(t)) = \nabla_{\bftheta} \hu(\bx;\bftheta(t))\cdot \dbftheta(t) - f(\bx,\hu(\cdot;\bftheta(t))),
\label{eq:Res}
\end{equation}
where $\dbftheta(t)$ is the result of applying the chain rule to the time derivative of $\hu(\cdot; \bftheta(t))$. Define $\mathcal{T}_{\bftheta(t)}$ to be the space spanned by the component functions of $\nabla_{\bftheta}\hu(\cdot; \bftheta(t))$, 
\[
\mathcal{T}_{\bftheta(t)} = \operatorname{span}\{\partial_{\theta_1(t)}\hat{u}(\cdot; \bftheta(t)), \dots, \partial_{\theta_N(t)}\hat{u}(\cdot; \bftheta(t))\}\,.
\]
The space $\mathcal{T}_{\bftheta(t)}$ is the tangent space of the manifold induced by the parameterization $\hu$ at parameter $\bftheta(t)$; we refer to \cite{Lubich2008} for details. We then seek $\dbftheta(t)$ such that it leads to a residual that is orthogonal in the following sense
\begin{equation}\label{eq:Prelim:GalerkinProj}
\langle \partial_{\theta_i(t)}\hat{u}(\cdot; \bftheta(t)), r_t(\cdot; \bftheta(t), \dbftheta(t))\rangle_{\nu} = 0\,, \qquad i = 1, \dots, \nh\,,
\end{equation}
where $\nu$ is a measure that is fully supported on $\Xcal$ and $\langle \cdot, \cdot \rangle_{\nu}$ denotes the corresponding $L^2$ inner product. After transformations, the system of equations \eqref{eq:Prelim:GalerkinProj} can be represented as
\begin{equation}
\bfM(\bftheta(t)) \dbftheta(t) = \bfF(\bftheta(t))\,,
\label{eq:MF_ODE}
\end{equation}
with the operators defined as
\[
[M(\bftheta(t))]_{ij} = \langle \nabla_{\theta_i(t)}\hu(\cdot; \bftheta(t)), \nabla_{\theta_j(t)}\hu(\cdot; \bftheta(t))\rangle_{\nu}\,,\qquad i,j = 1, \dots, \nh\,,
\]
and
\[
[F(\bftheta(t))]_i = \langle \nabla_{\theta_i(t)}\hu(\cdot; \bftheta(t)), f(\cdot; \hu(\cdot; \bftheta(t)))\rangle_{\nu}\,,\qquad i = 1, \dots, \nh\,.
\]
Thus, to obtain a numerical solution $\hu(\cdot; \bftheta(t))$ of the equation \eqref{eq:PDE}, the dynamical system \eqref{eq:MF_ODE} is integrated forward in time with a numerical time-integration scheme. 

\begin{figure}[t]
\begin{tabular}{cc}
\resizebox{0.50\columnwidth}{!}{\large\input{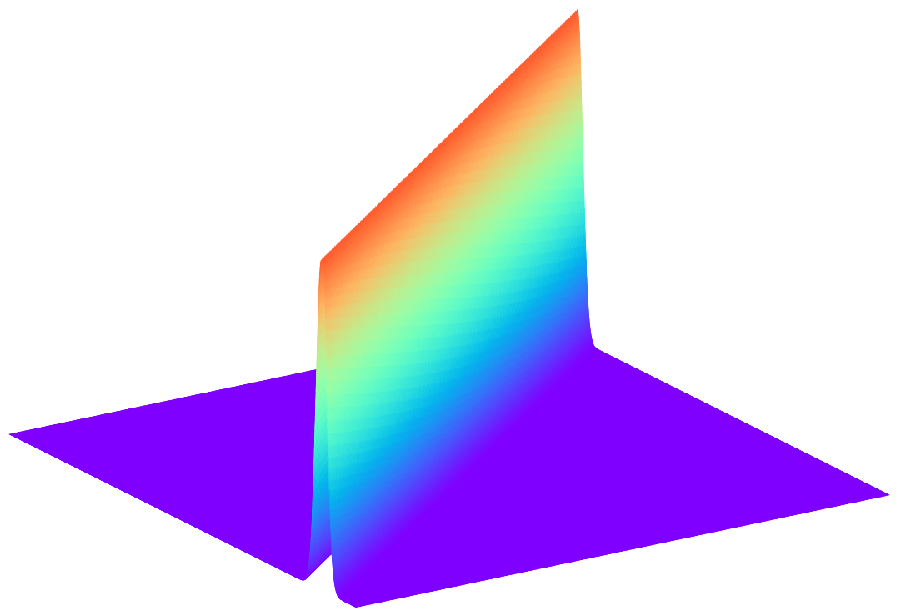}}~~~~~~~ & \resizebox{0.45\columnwidth}{!}{\Large\input{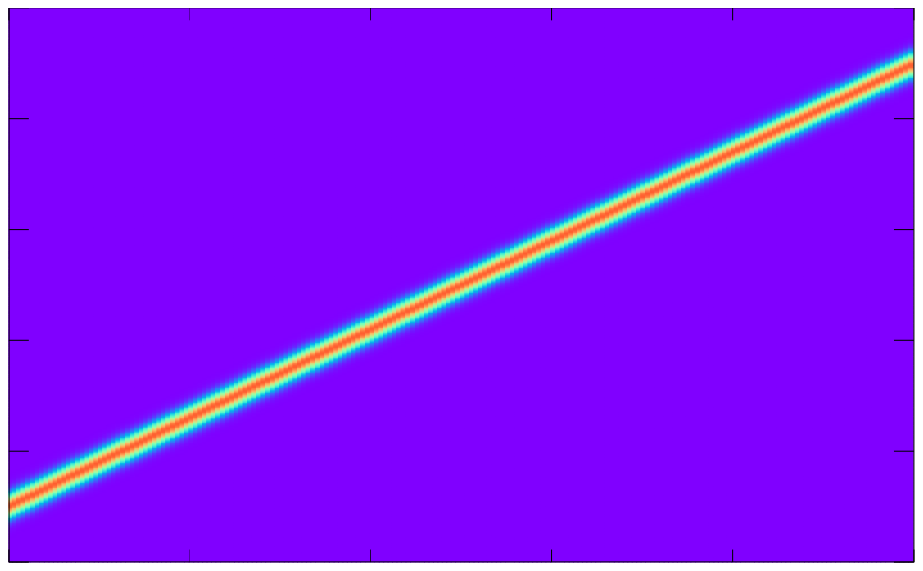}}\\
{\scriptsize (a) advecting wave} & {\scriptsize (b) time-space view of wave}
\end{tabular}
\caption{Features of solutions of transport-dominated problems are often local in the spatial domain. The spatial locality makes un-informed sampling such as uniform sampling inefficient for estimating the training loss of fitting nonlinear parametrizations, which is analogous to the challenge of estimating the operators $\bfM$ and $\bfF$ of the system \eqref{eq:MF_ODE} of our sequential-in-time approach.} 
\label{fig:Motivation}
\end{figure}

\subsection{Problem formulation}
\label{sec:Prelim:Prob}
For a limited class of PDEs and corresponding specific parametrizations, the operators $\bfM$ and $\bfF$ of the dynamical system shown in \eqref{eq:MF_ODE} can be derived analytically, which is demonstrated in \cite{Lubich2008,lasser_lubich_2020,doi:10.1137/21M1415972} in the context of time-dependent nonlinear parametrizations. In many cases, however, the operators $\bfM$ and $\bfF$ need to be estimated numerically. A common approach, especially if the spatial domain $\Xcal$ is higher dimensional, is to resort to Monte Carlo estimation from $\nm$ samples $\bx^{(1)}, \dots, \bx^{(\nm)} \in \Xcal$ of the distribution $\nu$.
However, Monte Carlo estimators based on un-informed sampling can have a high variance, for example if the solution $u$ has local features in the spatial domain and these local features are transported over time. Notice that this is a problem where nonlinear parametrizations can achieve a drastically faster error decay than linear parametrizations as discussed in the survey in \cite{P22AMS} and the introduction in Section~\ref{sec:Intro}. To see the sampling issue for such a problem, consider Figure~\ref{fig:Motivation}a that shows an advecting wave governed by the linear advection equation in one spatial dimension. By plotting the wave dynamics in the time-space domain in Figure~\ref{fig:Motivation}b, it can be seen that an un-informed sampling such as a uniform sampling in the time-space domain is inefficient in the sense that only few samples will be in regions where the solution is non-zero. While this is a toy example that is meant to demonstrate the challenge of estimating the operators of the system \eqref{eq:MF_ODE}, it is representative for a wide range of transport-dominated problems that have local features in the spatial domain. We stress that similar observations hold for global time-space collocation approaches, where an un-informed sampling can require a large number of samples to accurately estimate the residual norm over the time-space domain.

\begin{figure}
\resizebox{1.0\columnwidth}{!}{\fbox{\small\begin{tikzpicture}
\node (labelParameter) {parameters:};
\node[below = of labelParameter.east,left] (labelParticles) {particles:};

\node[right of = labelParameter, node distance=1in] (theta0) {$\bftheta_{t_0}$};
\node[below of = theta0] (x0) {$\mu_0 \sim \{\bfx_{t_0}^{(i)}\}_{i = 1}^m$};

\node[right of = theta0, node distance=1.5in] (theta1) {$\bftheta_{t_1}$};
\draw[->] (theta0.east) -- (theta1.west);

\node[right of = x0, node distance=1.5in] (x1) {$\{\bfx_{t_1}^{(i)}\}_{i = 1}^m$};
\draw[->] (x0.east) -- (x1.west);

\node[right of = theta1, node distance=1.5in] (theta2) {$\bftheta_{t_2}$};
\draw[->] (theta1.east) -- (theta2.west);
\node[right of = x1, node distance=1.5in] (x2) {$\{\bfx_{t_2}^{(i)}\}_{i = 1}^m$};
\draw[->] (x1.east) -- (x2.west);

\node[right of = theta2, node distance=1.5in] (theta3) {$\dots$};
\draw[->] (theta2.east) -- (theta3.west);
\node[right of = x2, node distance=1.5in] (x3) {$\dots$};
\draw[->] (x2.east) -- (x3.west);
\end{tikzpicture}}}\vspace*{-1.29cm}
\hspace*{-3.6cm}~\resizebox{1.4\columnwidth}{!}{\input{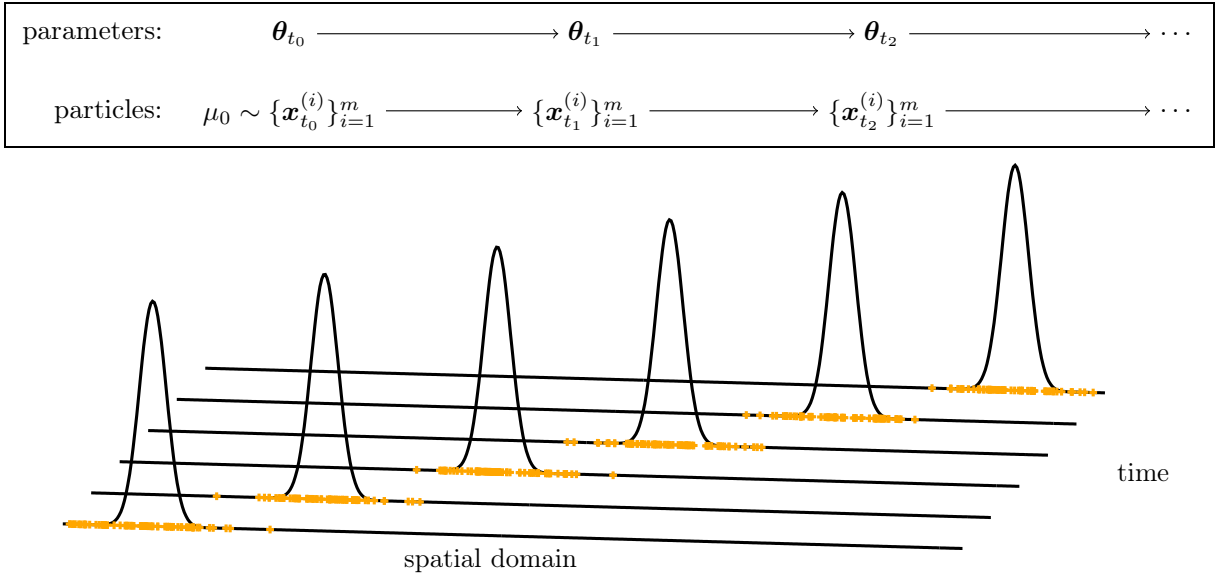}}
\vspace*{-1.29cm}\caption{The proposed Neural Galerkin schemes with dynamic particles couple the parameters $\bftheta(t)$ given by the Dirac-Frenkel variational principle with the ensemble of particles $\{\bfx_t^{(i)}\}_{i = 1}^{\nm}$ so that few particles are sufficient for accurately estimating the operators $\bfM_t$ and $\bfF_t$ of the system \eqref{eq:ANG:ParameterDyn}, even for problems with high-dimensional spatial domains and when there are local features that travel over time such as waves and other coherent structures.} 
\label{fig:Overview}
\end{figure}

\section{Neural Galerkin schemes with coupled parameter and particle dynamics}
\label{sec:adaptive_samp}
We introduce Neural Galerkin schemes with dynamic particles. The particles empirically represent time-dependent measures that are coupled to the parameters $\bftheta(t)$ that represent numerical approximations of the PDE solution fields. 
A coupled dynamical system integrates forward in time particles and parameters together such that few particles are sufficient for accurately estimating the inner products with the residual in the Dirac-Frenkel variational problem \eqref{eq:Prelim:GalerkinProj} and the parameters provide an approximation of the solution field in the Dirac-Frenkel variational sense; see Figure~\ref{fig:Overview}.  

This section is structured as follows. Section~\ref{sec:adaptive_IS} and Section~\ref{sec:ANG:Select} introduce Neural Galerkin schemes with projections that rely on time-dependent measures. The dynamics for the particles that approximate the time-dependent measure are derived in Section~\ref{sec:ANG:CoupledSys}, which then leads to the computational procedure of Neural Galerkin schemes with dynamic particles in Section~\ref{sec:ANG:Couple}. The computational procedure is summarized in algorithmic form in Section~\ref{sec:ANG:Comp}.

\subsection{Projections with time-dependent inner products}
\label{sec:adaptive_IS}
Instead of projecting with respect to the inner product corresponding to a fixed measure $\nu$ as in \eqref{eq:Prelim:GalerkinProj}, we introduce a time-dependent inner product $\langle \cdot, \cdot \rangle_{\mu_t}$ that is adapted over time with the time-dependent measure $\mu_t$. 
The projection with respect to $\mu_t$ is
\begin{equation}\label{eq:TNG:TimeGalerkinProj}
\langle \partial_{\theta_i(t)}\hat{u}(\cdot; \bftheta(t)), r_t(\cdot; \bftheta(t), \dbftheta(t))\rangle_{\mu_t} = 0\,, \qquad i = 1, \dots, \nh\,,
\end{equation}
with the corresponding $\bfM_t(\bftheta(t))$ and $\bfF_t(\bftheta(t))$ at time $t$  given by 
\[
[M_{t}(\bftheta(t))]_{i,j} = \langle \nabla_{\theta_i(t)}\hu(\cdot; \bftheta(t)), \nabla_{\theta_j(t)}\hu(\cdot; \bftheta(t))\rangle_{\mu_t}\,,\qquad i,j = 1, \dots, \nh\,,
\]
and
\[
[F_{t}(\bftheta(t))]_i = \langle \nabla_{\theta_i(t)}\hu(\cdot; \bftheta(t)), f(\cdot; \hu(\cdot; \bftheta(t)))\rangle_{\mu_t}\,,\qquad i = 1, \dots, \nh\,.
\]
Notice that $\bfM_t$ and $\bfF_t$ now depend on time and lead to the Neural Galerkin system
\begin{equation}
\bfM_t(\bftheta(t))\dbftheta(t) = \bfF_t(\bftheta(t))\,.
\label{eq:ANG:ParameterDyn}
\end{equation}
The aim is to adapt the measure $\mu_t$ over time $t$ so that the inner product \eqref{eq:TNG:TimeGalerkinProj} can be numerically estimated efficiently. If the nonlinear parametrization $\hu(\cdot; \bftheta(t))$ is so expressive that the component functions $\partial_{\theta_1(t)}\hu(\cdot; \bftheta(t))$, $\dots$, $\partial_{\theta_{\nh}(t)}\hu(\cdot; \bftheta(t))$ of the gradient of $\hu(\cdot; \bftheta(t))$ span an $\nh$-dimensional space and there exists $\bftheta(t)$ and $\dbftheta(t)$ with norm $|\langle r_t(\cdot; \bftheta(t), \dbftheta(t)), r_t(\cdot; \bftheta(t), \dbftheta(t)) \rangle_{\mu_t}| = \|r_t(\cdot; \bftheta(t), \dbftheta(t))\|_{\mu_t}^2 = 0$, then the corresponding $\bftheta(t)$ and $\dbftheta(t)$ will lead to a zero residual norm with respect to any other measure that has full support over $\Xcal$. Furthermore, under the expressiveness condition stated above, solving \eqref{eq:TNG:TimeGalerkinProj} for any measure $\mu_t$ implies that the residual norm is zero for any other measure.   
Thus, under these idealized conditions, we are free to choose a measure $\mu_t$ that is computationally convenient. In the following, we will select $\mu_t$ such that the inner product \eqref{eq:TNG:TimeGalerkinProj} can be accurately estimated with few samples. Notice that even computing the inner product in \eqref{eq:TNG:TimeGalerkinProj} with an approximation of $\mu_t$ is sufficient to obtain a zero residual norm with respect to the actual $\mu_t$. In numerical practice, the residual norm is typically not exactly zero; however, as long as the residual norm is kept small, similar arguments hold at least heuristically and provide a motivation for the following.

\subsection{Selecting the time-dependent measure}\label{sec:ANG:Select}
The Neural Galerkin schemes proposed in this work adapt the space $\mathcal{T}_{\bftheta(t)}$ over time following the Dirac-Frenkel variational principle as discussed in Section~\ref{sec:Prelim:NG} and additionally adapt the inner product $\langle \cdot, \cdot \rangle_{\mu_t}$ with respect to which the residual is set to zero. 

To this end, we will consider schemes that use  an inner-product with respect to evolving measures $\mu_t$ which track  Gibbs measures of the form 
\begin{equation}\label{eq:ANG:GibbsMeasure}
\mu^G_t(\mathrm d\bfx) = Z_{\bftheta(t), \dbftheta(t)}^{-1} \exp\left(-V_{\bftheta(t), \dbftheta(t)}(\bfx)\right) \mathrm d\bfx,\qquad Z_{\bftheta(t), \dbftheta(t)} = \int_{\mathcal{X}} \exp\left(-V_{\bftheta(t), \dbftheta(t)}(\bfx)\right) \mathrm d\bfx\,,
\end{equation}
where the potential $V_{\bftheta(t), \dbftheta(t)}$ depends on the parameter $\bftheta(t)$ and its time derivative $\dbftheta(t)$. These Gibbs measures are general enough for our purpose.

There are multiple options for choosing $\mu^G_t$. One option is setting 
\begin{equation}\label{eq:PNG:ResMu}
\mu^G_t \propto |r_t(\cdot; \bftheta(t), \dbftheta(t))|^{2\gamma}\nu(\cdot)^{\gamma}\,,
\end{equation}
where $\gamma > 0$ is a tempering parameter.
By construction, the samples from the measure $\mu^G_t$ defined in~\eqref{eq:PNG:ResMu}  are concentrated where the squared residual is high and thus typically few samples are sufficient to accurately estimate the inner products for the projections in \eqref{eq:TNG:TimeGalerkinProj}.

Another option is 
\begin{equation}
\mu^G_t \propto |\hu(\cdot; \bftheta(t))|^{\gamma}\nu(\cdot)^{\gamma}\,,
\label{eq:NGDyn:PropMagnitudeSolution}
\end{equation}
with tempering parameter $\gamma > 0$, which means that samples are drawn in regions of the spatial domain where the numerical solution $\hu$ has large values. We will later consider the Fokker-Planck equation where the solution $\hu$ is approximating a probability density function and thus sampling proportional to it leads to efficient estimators of the Neural Galerkin operators $\bfM_t$ and $\bfF_t$ of \eqref{eq:ANG:ParameterDyn}.

\subsection{Coupled system for the concurrent evolution of the parameters and the time-dependent measure}\label{sec:ANG:CoupledSys}

There are several possibilities to concurrently evolve an approximation $\mu_t$ of the measure $\mu^G_t$ and parameters~$\bftheta(t)$ corresponding to the solution field evolve.

\subsubsection{Langevin dynamics}
A first option is to evolve the parameters $\bftheta(t)$ and measure $\mu_t$ concurrently via the system of equations
\begin{equation}
\begin{aligned}
\bfM_t(\bftheta(t))\dbftheta(t) & = \bfF_t(\bftheta(t))\,,\\
\partial_t \mu_t & = \alpha \nabla \cdot (\nabla\mu_t + \mu_t\nabla V_{\bftheta(t), \dbftheta(t)})\,,
\end{aligned}\label{eq:ANG:CoupledSystem}
\end{equation}
where $\alpha>0$ is a parameter that controls the separation of timescale between the evolution of the parameters $\bftheta(t)$ and the measure $\mu_t$. In the limit as $\alpha\to\infty$, the measure evolves infinitely faster than the parameters, and since the equation for $\mu_t$ in~\eqref{eq:ANG:CoupledSystem} is the gradient flow in Wasserstein-2 metric over the Kullback-Leibler divergence of $\mu_t$ from the Gibbs measure associated with $V_{\bftheta(t), \dbftheta(t)}$, it realizes~\eqref{eq:ANG:GibbsMeasure} at all times. However, it is important to note that \eqref{eq:ANG:CoupledSystem} leads to zero-residual dynamics for any $\alpha>0$, as long as the nonlinear parametrization of the solution is expressive enough, as discussed at the end of Section~\ref{sec:adaptive_IS}.  In practice, \eqref{eq:ANG:CoupledSystem} can be implemented by  approximating the distribution $\mu_t$ by its empirical distribution over $\nm$ particles $\bfx^{(1)}(t), \dots, \bfx^{(\nm)}(t)$ as
\begin{equation}\label{eq:ANG:EmpDist}
\hat{\mu}_t = \frac{1}{\nm}\sum_{i = 1}^\nm \delta_{\bfx^{(i)}(t)}\,.
\end{equation}
This leads to the coupled system
\begin{equation}
\begin{aligned}
\hat{\bfM}_t(\bftheta(t))\dbftheta(t) & = \hat{\bfF}_t(\bftheta(t))\,,\\
\mathrm d\bfx^{(i)}(t) &= -\alpha\nabla V_{\bftheta(t), \dbftheta(t)}(\bfx^{(i)}(t))\mathrm d t + \sqrt{2\alpha}dW^{(i)}(t)\,,\qquad i = 1, \dots, \nm\,,
\end{aligned}\label{eq:ANG:CoupledSystemDisc}
\end{equation}
where $\{W^{(i)}(t)\}_{i=1}^m$ are $m$ independent Wiener processes in $\mathbb{R}^d$ and the particles $\bfx^{(1)}(t), \dots, \bfx^{(m)}(t)$ are used to estimate $\bfM_t$ and $\bfF_t$ as
\begin{equation}
\begin{aligned}
[\hat{M}_{t}(\bftheta(t))]_{l,j} = & \frac1m \sum_{i=1}^m \nabla_{\theta_l(t)}\hu(\bfx^{(i)}(t); \bftheta(t)) \nabla_{\theta_j(t)}\hu(\bfx^{(i)}(t); \bftheta(t))\,,\qquad l,j = 1, \dots, \nh\,,\\
[\hat{F}_{t}(\bftheta(t))]_j = &\frac1m \sum_{i=1}^m \nabla_{\theta_j(t)}\hu(\bfx^{(i)}(t); \bftheta(t)) f(\bfx^{(i)}(t); \hu(\cdot; \bftheta(t)))\,,\qquad j = 1, \dots, \nh\,.
\end{aligned}
\label{eq:ANG:EstimatedMF}
\end{equation}

\subsubsection{Stein variational gradient descent}\label{sec:ANG:SVGD}
We can also use Stein variational gradient descent (SVGD) \cite{NIPS2016_b3ba8f1b} to derive dynamical systems for the particles. The SVGD method deterministically approximates the gradient flow \eqref{eq:ANG:CoupledSystem} in a reproducing kernel Hilbert space with kernel $\mathcal{K}: \mathbb{R}^d \times \mathbb{R}^d \to \mathbb{R}$. Approximating the gradient flow in an RKHS can be beneficial when the gradients of the potential $V$ are noisy, which is typically the case when taking the gradient of the potential involves differentiating a neural-network approximation. 

For using SVGD, we replace the equation for the measure $\mu_t$ in  \eqref{eq:ANG:CoupledSystem}  with
\begin{equation}\label{eq:ANG:SVGDPDE:0}
\partial_{t}\mu_t = \alpha \nabla \cdot\left(\mu_t \mathbb{E}_{\bfx' \sim \mu_t}[\mathcal{K}(\bfx', \bfx)\nabla V_{\bftheta(t), \dbftheta(t)}(\bfx') - \nabla_1\mathcal{K}(\bfx', \bfx)]\right)\,,
\end{equation}
where we denote with $\nabla_1\mathcal{K}$ the gradient of $\mathcal{K}$ in its first argument. 
In practice, we approximate $\mu_t$ with its empirical distribution as in \eqref{eq:ANG:EmpDist} and obtain from \eqref{eq:ANG:SVGDPDE:0} the particle dynamics 
\begin{equation}
\label{eq:SVGDParticleFlow:b:0}
\frac{\mathrm d}{\mathrm dt}\bfx^{(i)}(t) = \frac\alpha m \sum_{l=1}^m [\nabla_1 \mathcal{K}(\bfx^{(l)}(t), \bfx^{(i)}(t)) - \mathcal{K}(\bfx^{(l)}(t), \bfx^{(i)}(t))\nabla V_{\bftheta(t), \dbftheta(t)}(\bfx^{(l)}(t))]\,,\quad i = 1, \dots, \nm\,.
\end{equation}
The choice of the kernel $\mathcal{K}$ plays an important role in SVGD in theory and practice. In particular, a poor choice of the kernel can prevent convergence of \eqref{eq:ANG:SVGDPDE:0} to the Gibbs measure \eqref{eq:ANG:GibbsMeasure} with potential $V$ \cite{doi:10.1137/18M1187611} and even a failure of detecting non-convergence in practical settings \cite{pmlr-v70-gorham17a}. However, as stated in Section~\ref{sec:adaptive_IS}, if the nonlinear parametrization is expressive enough, it is sufficient in our case to compute the projections in \eqref{eq:TNG:TimeGalerkinProj} with respect to an approximation of the measure $\mu_t^G$.

\subsection{Scheme for concurrently evolving parameter and particle dynamics}\label{sec:ANG:Couple}
We now propose a numerical scheme to integrate the coupled system of parameter dynamics \eqref{eq:ANG:ParameterDyn} and particle dynamics given by, e.g., SVGD as in \eqref{eq:SVGDParticleFlow:b:0}. Motivated by the heterogeneous multiscale method \cite{cms/1118150402,E20095437,abdulle_weinan_engquist_vanden-eijnden_2012}, we propose to move forward the ensemble of particles $\{\bfx^{(i)}(t)\}_{i = 1}^{\nm}$ with a faster time scale than the parameters $\bftheta(t)$ over time $t$ given by the PDE \eqref{eq:PDE}. 

\subsubsection{Discretization in time}
We discretize the dynamical system \eqref{eq:ANG:ParameterDyn} of the parameters $\bftheta(t)$ in time. Let $\delta t_k > 0$ be the time-step size at time step $k = 0, 1, 2, ...$ so that at time step $k + 1$ the time is $t_{k + 1} = t_k + \delta t_k$. Correspondingly, let $\Delta \bftheta_k$ be the update that is applied at time step $k$ to obtain $\bftheta_{k + 1} = \bftheta_k + \delta t_k\Delta \bftheta_k$. The initial parameter is $\bftheta_0$, which is obtained by fitting $\hu(\cdot; \bftheta_0)$ to the initial condition $u_0$.

The update $\Delta \bftheta_k$ at time step $k$ is obtained by solving the time-discrete system
\begin{equation}
\hbfM_k(\bftheta_k, \Delta\bftheta_k) \Delta \bftheta_k = \hbfF_k(\bftheta_k, \Delta\bftheta_k)\,,
\label{eq:ANG:TimeDiscSystem}
\end{equation}
which is obtain after discretizing the time-continuous system \eqref{eq:ANG:ParameterDyn} and estimating the operators with Monte Carlo over the particles $\bfx_k^{(1)}, \dots, \bfx_k^{(m)}$ at time step $k$ analogous to \eqref{eq:ANG:EstimatedMF}. Notice that \eqref{eq:ANG:TimeDiscSystem} describes potentially nonlinear equations in $\Delta\bftheta_k$. For example, if the time-continuous system \eqref{eq:ANG:ParameterDyn} is discretized with implicit time integration schemes such as the backward Euler method, then $\hbfM_k$ and $\hbfF_k$ can depend nonlinearly on $\Delta\bftheta_k$. In contrast, for explicit schemes, system \eqref{eq:ANG:TimeDiscSystem} is linear in $\Delta\bftheta_k$; we refer to \cite{NG22} for more details.

To keep the notation manageable, we focus on the Runge-Kutta 4 (RK4) scheme \cite{DORMAND198019} here. The corresponding update is given by
\[
\Delta\bftheta_k = \frac{1}{6}\left(\Delta\bftheta_k^{(1)} + 2\Delta\bftheta_k^{(2)} + 2 \Delta\bftheta_k^{(3)} + \Delta\bftheta_k^{(4)}\right)\,,
\]
where
\begin{equation}
\begin{aligned}
\hat{\bfM}_k\left(\bftheta_k, \Delta\bftheta_k^{(1)}\right) \Delta\bftheta_k^{(1)} = &\hat{\bfF}_k\left(\bftheta_k, \Delta\bftheta_k^{(1)}\right)\,,\\
\hat{\bfM}_k\left(\bftheta_k + \frac{\delta t_k}{2}\Delta\bftheta_k^{(1)}, \Delta\bftheta_k^{(2)}\right) \Delta\bftheta_k^{(2)} = &\hat{\bfF}_k\left(\bftheta_k + \frac{\delta t_k}{2}\Delta\bftheta_k^{(1)}, \Delta\bftheta_k^{(2)}\right)\,,\\
\hat{\bfM}_k\left(\bftheta_k + \frac{\delta t_k}{2}\Delta\bftheta_k^{(2)}, \Delta\bftheta_k^{(3)}\right) \Delta\bftheta_k^{(3)} = & \hat{\bfF}_k\left(\bftheta_k + \frac{\delta t_k}{2}\Delta\bftheta_k^{(2)}, \Delta\bftheta_k^{(3)}\right)\,,\\
\hat{\bfM}_k\left(\bftheta_k + \delta t_k\Delta\bftheta_k^{(3)}, \Delta\bftheta_k^{(4)}\right) \Delta\bftheta_k^{(4)} = &\hat{\bfF}_k\left(\bftheta_k + \delta t_k\Delta\bftheta_k^{(3)}, \Delta\bftheta_k^{(4)}\right)\,.
\end{aligned}
\label{eq:ANG:RK4Equations}
\end{equation}
The estimated matrices are
\begin{equation}
\begin{aligned}
[\hat{M}_k(\bftheta, \Delta\bftheta)]_{l,j} &= \frac{1}{\nm}\sum_{i=1}^{\nm}\nabla_{\theta_l} \hu(\bx_k^{(i)}; \bftheta) \nabla_{\theta_i} \hu(\bx_k^{(i)}; \bftheta)\,, \qquad l,j = 1, \dots, \nh\,,\\
[\hat{F}_k(\bftheta, \Delta\bftheta)]_j &= \frac{1}{\nm}\sum_{i=1}^{\nm}\nabla_{\theta_j} \hu(\bx_k^{(i)}; \bftheta) f(\bx_k^{(i)}, \hu(\cdot; \bftheta))\,,\qquad j = 1, \dots, \nh\,,
\end{aligned}
\label{eq:ANG:EstMK}
\end{equation}
with the particles $\bfx_k^{(1)}, \dots, \bfx_k^{(m)}$ at time step $k$. Notice that $\Delta \bftheta$ only formally enters in \eqref{eq:ANG:EstMK} but does not change $\hbfM_k$ and $\hbfF_k$ and thus the four equations in \eqref{eq:ANG:RK4Equations} are linear in the updates $\Delta\bftheta_k^{(1)}, \dots, \Delta\bftheta_k^{(4)}$. This means that at each time step with RK4, four systems of linear equations are solved.

\subsubsection{Predictor-corrector scheme to update particles}
To obtain the particles $\bfx_k^{(1)}, \dots, \bfx_k^{(n)}$ for computing \eqref{eq:ANG:EstMK} at time step $k$, we first take a predictor step because the potential $V_{\bftheta_k, \Delta\bftheta_k}$ at time step $k$ depends on $\bftheta_k$ and $\Delta \bftheta_k$, which is unavailable. With a predictor step we compute $\Delta \bftheta_k^{(P)}$ as an approximation of $\Delta \bftheta_k$ to use in adapting the particles. The predictor system is obtained with a forward Euler discretization of \eqref{eq:ANG:ParameterDyn}, so that
\begin{equation}
\begin{aligned}
[\hat{M}_{k}^{(P)}]_{l,j} &=  \frac{1}{\nm}\sum_{i=1}^{\nm}\nabla_{\theta_l} \hu(\bx_{k-1}^{(i)}; \bftheta_k) \nabla_{\theta_j} \hu(\bx_{k-1}^{(i)}; \bftheta_k)\,,\qquad l,j = 1, \dots, \nh\,,\\ 
[\hat{F}_k^{(P)}]_j &= \frac{1}{\nm}\sum_{i=1}^{\nm}\nabla_{\theta_j} \hu(\bx_{k-1}^{(i)}; \bftheta_k) f(\bx_{k-1}^{(i)}, \hu(\cdot; \bftheta_k))\,,\qquad j = 1, \dots, \nh\,,
\end{aligned}
\label{eq:ANG:PredictorSys}
\end{equation}
which leads to a linear system of equations in $\Delta \bftheta_k^{(P)}$. Notice that the estimates $\hbfM_k^{(P)}$ and $\hbfF_k^{(P)}$ are based on the particles $\bfx_{k - 1}^{(1)}, \dots, \bfx_{k - 1}^{(\nm)}$ from time step $k - 1$. 

We use the result of the predictor $\Delta \bftheta_k^{(P)}$ to define the potential $V_{\bftheta_k, \Delta \bftheta_k^{(P)}}$, which is then used to update the particles as described in Section~\ref{sec:ANG:CoupledSys}. The initial distribution is the empirical measure $\hat{\mu}_{k - 1}$ from the previous time step $k - 1$. Integrating the particle dynamics forward in time with a time-step size $\delta \tau > 0$ until particle end time gives the particles at time step $k$,
\[
\bfx_k^{(1)}, \dots, \bfx_k^{(\nm)}\,.
\]
The time-step size $\delta \tau = \delta t/\alpha$ with which the particles are integrated in time is controlled by the parameter $\alpha$ of equations \eqref{eq:ANG:CoupledSystemDisc} and \eqref{eq:SVGDParticleFlow:b:0}, respectively.
The particles $\bfx_k^{(1)}, \dots, \bfx_k^{(\nm)}$ are then used to compute \eqref{eq:ANG:EstMK} and solve for $\Delta \bftheta_k$ to obtain $\bftheta_{k + 1}$.
This process is repeated for all time steps $k = 1, 2, 3, \dots$.

\subsection{Computational procedure}\label{sec:ANG:Comp}
We summarize the computational procedure of Neural Galerkin schemes with dynamic particles in Algorithm~\ref{alg:coupling_dynamics}. The procedure requires as inputs the parameter $\bftheta_0$ that is obtained by fitting $\hu(\cdot; \bftheta_0)$ to the initial condition $u_0$ of the PDE \eqref{eq:PDE}. Other inputs are the number of time steps $K$ and the time-step size $\delta t$. For simplicity, we use a fixed time step size $\delta t = \delta t_k$ for $k=1,\dots,K$. Additionally, an initial ensemble of particles $\{\bfx_0^{(i)}\}_{i = 1}^{\nm}$ is required, which can be obtained by sampling from the distribution that is proportional to the absolute initial condition $|u_0|$. 

The procedure consists of a loop over the time steps $k = 1, \dots, K$ corresponding to the physical time $t$ and a nested loop to integrate the particle dynamics over time $\tau$ as described in Section~\ref{sec:ANG:CoupledSys}. At each time step $k = 1, \dots, K$, the predictor step gives $\Delta\bftheta_k^{(P)}$, which is then used to compute the potential $V_{\bftheta_k, \Delta \bftheta_k^{(P)}}$ and its gradient for updating the particles. The particles are initialized at time step $k$ with the ensemble of particles $\{\bfx_{k-1}^{(i)}\}_{i = 1}^{\nm}$ from the previous time step $k - 1$ and then updated to the ensemble of particles $\{\bfx_k^{(i)}\}_{i = 1}^{\nm}$ at time $k$. The particles $\bfx_k^{(1)}, \dots, \bfx_k^{(\nm)}$ are then used to estimate $\hbfM_k$ and $\hbfF_k$ to assemble the system \eqref{eq:ANG:TimeDiscSystem}, which is then integrated forward for one time step to obtain the update $\Delta\bftheta_k$ and ultimately $\bftheta_{k + 1} = \bftheta_k + \Delta\bftheta_k$.

\begin{algorithm}[bt]
\caption{Coupling parameter and particle dynamics in Neural Galerkin schemes}
\label{alg:coupling_dynamics}
\renewcommand{\algorithmicrequire}{\textbf{Input:}}
\begin{algorithmic}
\Require $\{\bfx_{0}^{(i)}\}_{i=1}^m$, $\bftheta_0$, $K$, $\delta t$
\For{$k \gets 1$ to $K$}                    
    \State {$\Delta\bftheta_k^{(P)}$ $\gets$ {predictor step \eqref{eq:ANG:PredictorSys} with the operators estimated with ensemble $\{\bfx_{k-1}^{(i)}\}_{i = 1}^{m}$}}
    \State {$\{\bfx_{k}^{(i)}\}_{i = 1}^{m}$ $\gets$ initialize with $\{\bfx_{k-1}^{(i)}\}_{i = 1}^{m}$ from time step $k-1$ and then update as in Section~\ref{sec:ANG:CoupledSys}}
   \State $\hbfM_k, \hbfF_k$ $\gets$ Estimate the operators with ensemble $\{\bfx_k^{(i)}\}_{i = 1}^{\nm}$
    \State {$\Delta\bftheta_k$ $\gets$ take a time step of system \eqref{eq:ANG:TimeDiscSystem} with the estimated operators $\hbfM_k$ and $\hbfF_k$} 
    \State {$\bftheta_{k+1}$ $\gets$ $\bftheta_k + \delta t\Delta\bftheta_k$} update parameters
\EndFor
\end{algorithmic}
\end{algorithm}

\section{Numerical experiments}
\label{sec:numerical_results}
In this section, we demonstrate Neural Galerkin schemes with coupled particles dynamics on three numerical examples. We start with the one-dimensional Korteweg-de Vries (KdV) equation, for which we plot the particle dynamics for demonstration purposes. We then consider transport equations and the Fokker-Planck equation in moderately high dimensions. We will demonstrate the effectiveness of our adaptive sampling scheme with dynamic particles by comparing it to sampling from a measure $\nu$ that is fixed over time. In all three examples, $\nu$ is chosen to be the uniform measure over the spatial domain $\Xcal$.

\subsection{Korteweg-de Vries equation (KdV)}\label{sec:NumExp:KdV}
Consider the KdV equation 
\begin{equation}
\partial_tu(t, x) + \partial_x^3u(t, x) + 6u(t, x)\partial_xu(t, x) = 0\,,\qquad (t, x) \in [0, 6] \times \mathcal{X}\,,
\label{eq:NumExp:KdV}
\end{equation}
over the spatial domain $\mathcal{X} = [-20, 40)$ and time range $t\in[0,6]$. We consider the same initial condition $u_0$ as in \cite{Taha:1984ka} for which an analytic solution to \eqref{eq:NumExp:KdV} is available. The solution consists of two solitons that approach each other, merge, and then separate again. We approximate the solution by imposing Dirichlet boundary conditions on \eqref{eq:NumExp:KdV} via a penalty residual on the left $x_l=-20$ and right $x_r=40$ boundary point. Consider the residual at the boundary
\begin{equation} r_t^\partial(\cdot; \bftheta(t),\dbftheta(t)) = \nabla_{\bftheta} \hu(\cdot;\bftheta(t))\cdot \dbftheta(t) - g(\cdot;\bftheta(t))\,,
    \label{eq:Res_bc}
\end{equation}
where the function $g \equiv 0$ is constant zero so that we penalize any deviation over time from the value of initial condition at the boundary points. Then the residual that we project following \eqref{eq:TNG:TimeGalerkinProj} is
\[
\bar{r}_t(\cdot; \bftheta(t), \dbftheta(t)) = r_t(\cdot; \bftheta(t), \dbftheta(t)) + \zeta(r_t^{\partial}(x_l; \bftheta(t), \dbftheta(t)) + r_t^{\partial}(x_r; \bftheta(t), \dbftheta(t)))
\]
where $\zeta > 0$ is an adjustable weight we put on the boundary condition penalty. In this experiment, we set it to $\zeta = 10^4$ and note that in later examples we consider a different way of imposing boundary conditions that does not require setting an adjustable weight.
For the nonlinear parameterization $\hu$, we use a fully connected feed-forward network with two hidden layers, five nodes per layer, and sigmoid activation functions. The total number of parameters is $\nh = 45$. We use the RK4 time integration scheme with a fixed time-step size $\delta t=10^{-4}$.

We compare our adaptive sampling scheme with dynamic particles to a static sampling that keeps the distribution fixed as the uniform distribution over the spatial domain $\mathcal{X}$. Our target measure is set to $\mu_t\propto |r_t(\cdot; \bftheta(t), \dbftheta(t))|^{2}\nu(\cdot)$ and the number of particles is $\nm = 100$, which is the same number of samples drawn from the uniform measure at each time $t$. Particles dynamics are imposed via SVGD as described in Section~\ref{sec:ANG:SVGD}. The kernel is the Gaussian kernel with bandwidth $0.05$, the step size for SVGD is $0.05$, and the number of SVGD steps is 500 at each physical time $t$. The tempering parameter is set to $\gamma=0.25$. SVGD particles are initialized proportional to the fitted initial solution $\hu(\cdot;\bftheta_0)$ at $t=0$. 

\begin{figure}
\centering
\begin{tabular}{cc}
    \resizebox{.5\linewidth}{!}{\Large \input{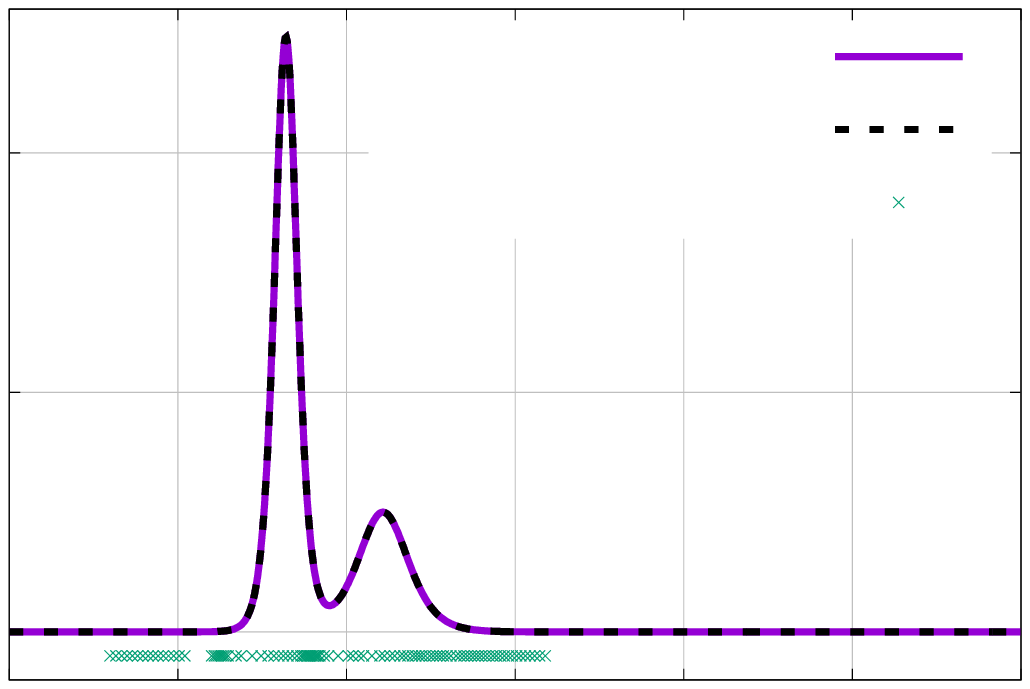}} &  \resizebox{.5\linewidth}{!}{\Large \input{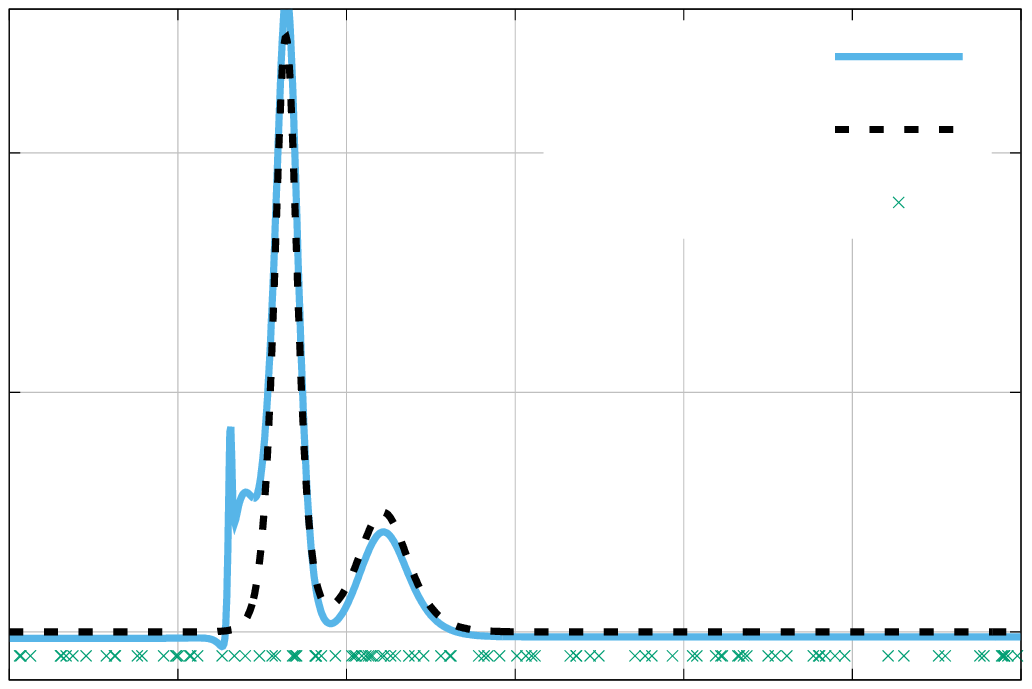}}\\
   \scriptsize (a) Neural Galerkin, time $t = 0.3$ & \scriptsize (b) static, time $t = 0.3$\\
    \resizebox{.5\linewidth}{!}{\Large \input{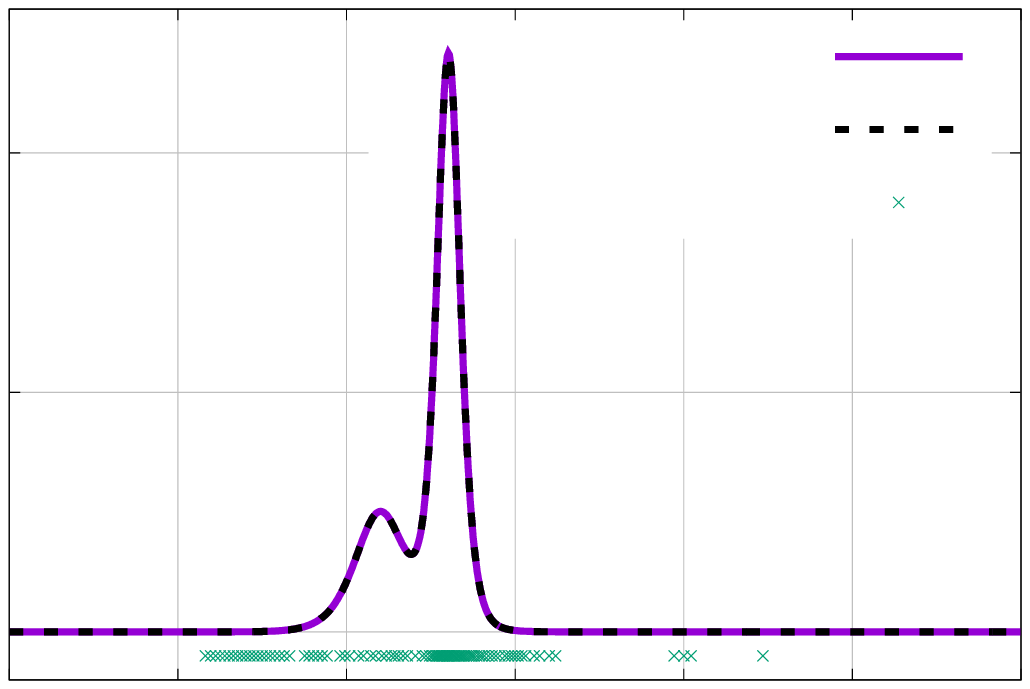}} &  \resizebox{.5\linewidth}{!}{\Large \input{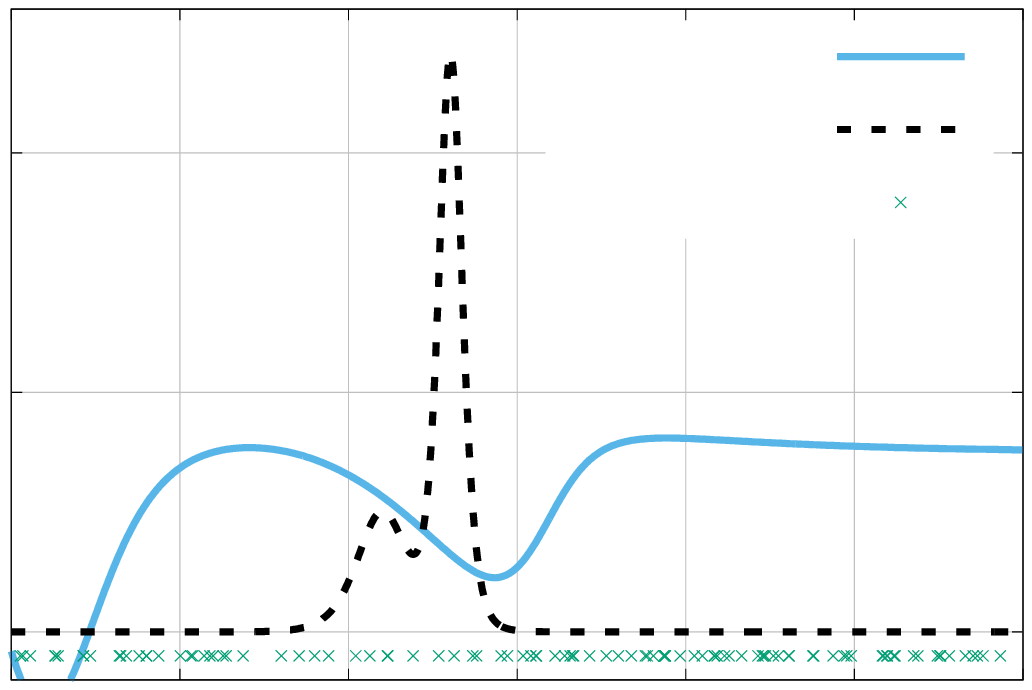}}\\
    \scriptsize (c) Neural Galerkin, time $t = 2.0$ & \scriptsize (d) static, time $t = 2.0$ \\
    \resizebox{.5\linewidth}{!}{\Large \input{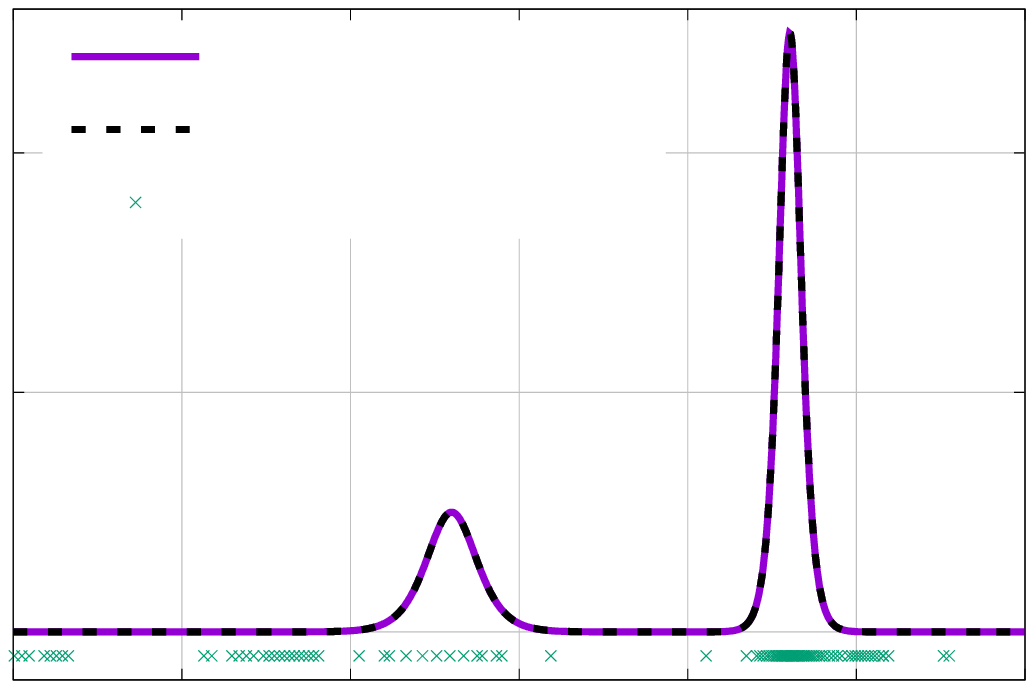}} &  \resizebox{.5\linewidth}{!}{\Large \input{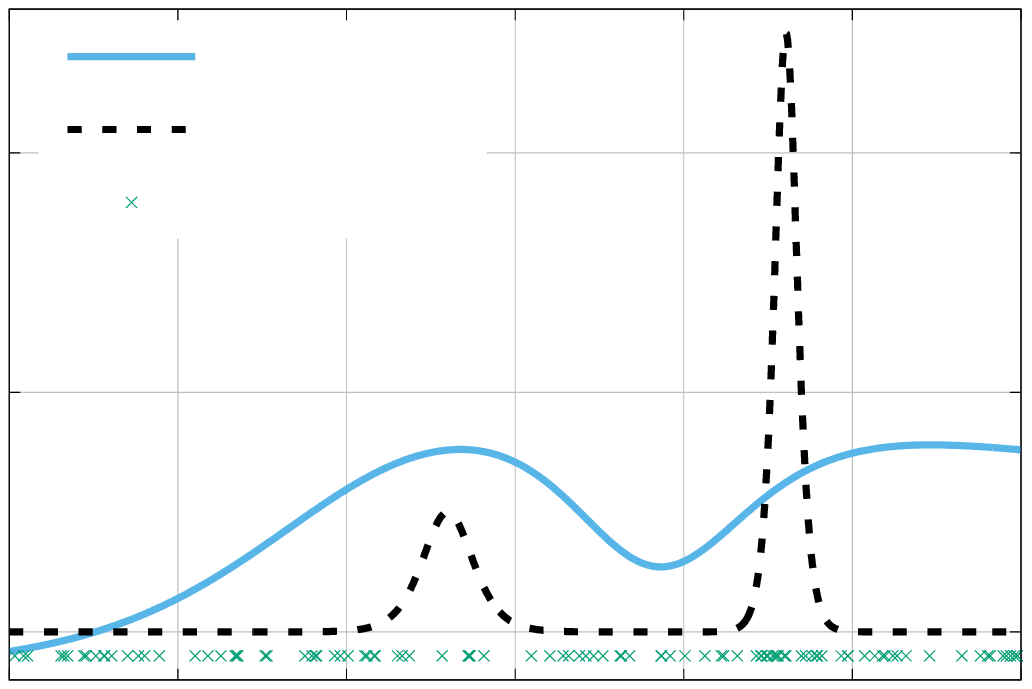}}\\
    \scriptsize (e) Neural Galerkin, time $t = 6.0$ & \scriptsize (f) static, time $t = 6.0$ 
\end{tabular}
    \caption{KdV: Neural Galerkin schemes with dynamic particles (left column) achieves an accurate approximation of the solution field with only $\nm = 100$ particles, whereas static sampling (right column) with $\nm = 100$ samples leads to inaccurate approximations after only a few time steps.} 
    \label{fig:KdV_adaptiveVSstatic}
\end{figure}

\begin{SCfigure}[50][ht]
\resizebox{.5\linewidth}{!}{\Large \input{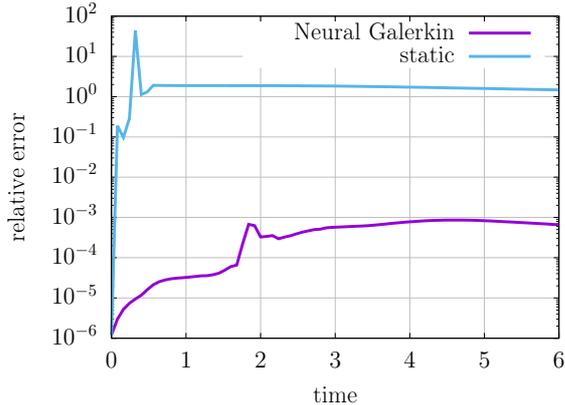}}
\caption{KdV: Neural Galerkin schemes with dynamic particles lead to orders of magnitude lower errors than approximations based on static sampling in this example.}
\label{fig:KdV_error}
\end{SCfigure}

Figure~\ref{fig:KdV_adaptiveVSstatic} compares the approximation of the solution obtained with Neural Galerkin schemes with dynamic particles to the approximation obtained with static sampling. Dynamically moving the particles leads to a good approximation of the solution that captures the local features of this problem. In contrast, the accuracy of the approximation based on uniform samples quickly deteriorates. A quantitative comparison between Neural Galerkin schemes with dynamic particles versus static approaches is shown in Figure~\ref{fig:KdV_error}. The relative $L^2$ error is computed with respect to the analytic solution over the spatial domain $\Xcal$. As shown in Figure~\ref{fig:KdV_error}, Neural Galerkin schemes with dynamic particles lead to orders of magnitude lower errors than static sampling.

\sidecaptionvpos{figure}{c}

\subsection{High-dimensional transport equations}
Consider the advection equation 
\[
\partial_tu(t, \bfx) + \ba(t)\cdot \nabla_\bx u(t, \bfx) = 0,\quad\quad u(0,\bfx) = u_0(\bfx)\,,
\]
over the five-dimensional spatial domain $\Xcal = [0, 10]^5$ for time range $t\in[0,1.2]$. The time-dependent transport coefficient is 
\begin{equation}
\ba(t) = \bc\odot(\sin(\pi t\ba_d) + 5/4)\,,
\end{equation}
where $\bc = [1,2,...,d]^T$, $\ba_d = 2+\frac{2}{d}[0,1,...,d-1]^T$, and $\odot$ denotes element-wise multiplication. The initial condition $u_0$ is a mixture of two non-isotropic Gaussian waves with means
\[
\bfmu_1 = \frac{11}{10}\begin{bmatrix}
    1\\ 1\\ \vdots\\ 1
\end{bmatrix},\quad\quad \bfmu_2 = \frac{3}{4}\begin{bmatrix}
    1.5-(-1)^1/(d+1)\\
    1.5-(-1)^2/(d+1)\\
    \vdots\\
    1.5-(-1)^d/(d+1)
\end{bmatrix}\,,
\]
and covariance matrices
\[
\bfSigma_1 = \frac{1}{200}\begin{bmatrix}
    2\\ & 4\\ & &\ddots\\ & & & 2d
\end{bmatrix},\quad\quad \bfSigma_2 = \frac{1}{200}\begin{bmatrix}
    d+1\\ & d\\ & & \ddots\\ & & & 2
\end{bmatrix}.
\]
We enforce that the solution is zero
at the origin with a penalty term that is applied similarly as described in Section~\ref{sec:NumExp:KdV}. The coefficient of the penalty term is $\zeta = 10^2$ in this example. The analytic solution of this equation can be derived via the method of characteristics; see \cite{NG22} for details. The analytic solution will serve as benchmark in this example. 
For the nonlinear parameterization, we use a fully connected feed-forward network with two hidden layers, fifteen nodes per layer, and sigmoid activation functions. The number of parameters is $\nh = 345$. We use the RK4 time integration scheme with a fixed time-step size $\delta t=10^{-3}$.

We compare Neural Galerkin with dynamic particles to uniform sampling with $m = 2500$ samples over the five-dimensional domain $\Xcal$. The target measure $\mu_t$ is proportional to the squared residual as in Section~\ref{sec:NumExp:KdV}. In this example, we take 300 SVGD steps at each time $t$ and set the kernel bandwidth and the SVGD step size to 0.1. The tempering parameter is set to $\gamma=0.25$.

\begin{figure}
\centering
    \includegraphics[width=\textwidth]{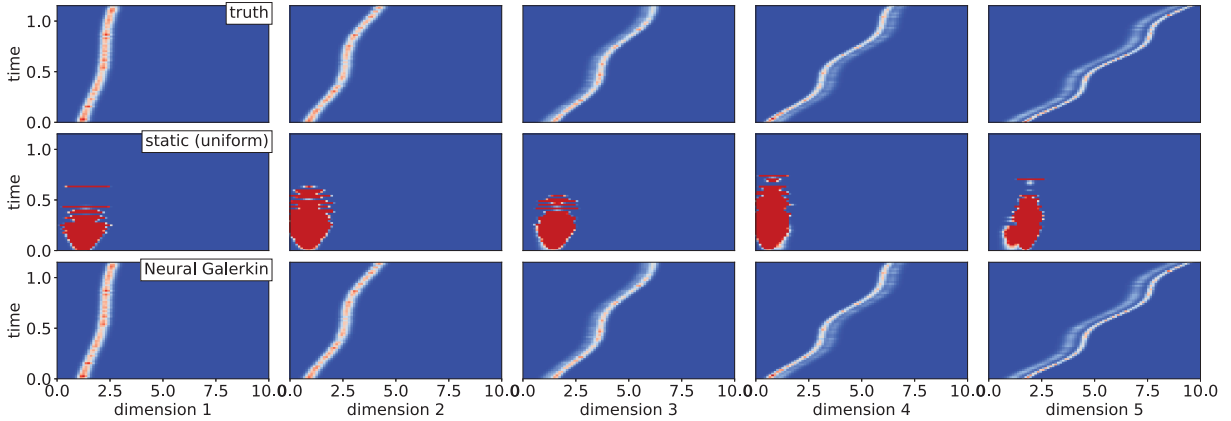}
    \caption{High-dimensional advection: This figure shows the marginals of the 
analytic solution (truth) and of the numerical solutions obtained with static sampling (uniform) and dynamic particles (Neural Galerkin). With only $\nm = 2500$ samples, Neural Galerkin with dynamic particles can predict well the local features of the solution in this moderately high-dimensional example. In contrast, static sampling based on the uniform distribution fails to provide meaningful predictions.}
    \label{fig:LinAdv_marginals}
\end{figure}

Figure~\ref{fig:LinAdv_marginals} compares the analytic solution to the approximate solutions obtained with dynamic particles and with static sampling and $\nm = 2500$ samples. The figure shows the marginals 
\[
u^{\text{marg}}_i(t, x) = \int_0^{10} \cdots \int_0^{10} u(t, x_1, \dots, x_{i - 1}, x, x_{i + 1}, \dots, x_d)\mathrm{d}x_1 \dots \mathrm d x_{i - 1} \mathrm d x_{i + 1} \dots \mathrm dx_d
\]
for dimensions $i = 1, \dots, d$ of the analytic solution, the numerical solution obtained with static sampling, and the numerical solution obtained with Neural Galerkin and dynamic particles. 
The results shown in Figure~\ref{fig:LinAdv_marginals} indicate that 2500 uniform samples over the five-dimensional spatial domain $\Xcal$ are insufficient to capture the local features and so lead to a poor approximation of the analytic solution. In contrast, if we dynamically adapt the particles over time, we obtain an approximation that is in close agreement with the analytic solution in this example.

\subsection{Non-Gaussian interaction systems described by the Fokker-Planck equation}
\label{sec:NumExp:Fokker_Planck}
We consider a system of $d$ interacting physical particles in a bounded domain $\mathcal{X}=[-3,11]^d\subset \mathbb{R}^d$ with positions $X_1(t),...,X_d(t)$ described by the stochastic differential equation (SDE)
\begin{equation}
    \mathrm dX_i = f_g(t,X_i)\mathrm dt + \sum_{j=1}^df_K(X_i,X_j)\mathrm dt + \sqrt{2D}\mathrm dW_i,\quad \textnormal{for } i=1,...,d,
    \label{eq:PT_particles}
\end{equation}
with one-body force $f_g:[0,\infty)\times\mathbb{R}\rightarrow \mathbb{R}$, a pairwise interaction term $f_K:\mathbb{R}\times\mathbb{R}\rightarrow\mathbb{R}$, the diffusion coefficient $D>0$, and independent Wiener processes $W_i$. 

\paragraph{Fokker-Planck equation} The joint density of the positions $X_1(t),...,X_d(t)$ is governed by the Fokker-Planck equation with homogeneous Dirichlet boundary conditions,
\begin{equation}
    \partial_tu(t, \bfx) = \sum_{i=1}^d -\partial_{x_i}(u(t, \bfx) h_i(t,\bx)) + D\partial^2_{x_i}u(t, \bfx)\,,
\end{equation}
with $h_i(t,\bx) = h_i(t,x_1,...,x_d) = f_g(t,x_i) + \sum_{j=1}^df_K(x_i,x_j)$. We consider this problem in $d=8$ dimensions, the one-body force $f_g(t,\bfx) = \frac{5 \times 10^{1/3}}{4}(\sin(\pi t)+\frac{3}{2})-\bfx$, interacting term $f_K(\bfx,\bfy) = \frac{1}{2d}(\bfy-\bfx)$, and diffusion coefficient $D=0.5$. The initial condition $u_0$ is an isotropic Gaussian density with mean $\frac{29}{10}+\frac{21}{10(d-1)}[0, 1, \dots, d-1]^T$ and variance $\sigma^2=0.1$.

\paragraph{Parametrization and enforcing boundary conditions} Since the solution $u$ represents a probability density function, we parameterize the potential $\log u$ of $u$, instead of $u$ directly. Additionally, we enforce homogeneous Dirichlet boundary conditions via a product structure so that we obtain the parametrization
\[
\hu(\bfx; \bftheta(t)) = \hu_{\text{BC}}(\bfx)\exp(\hu_{\text{potential}}(\bfx;\bftheta(t)))\,,
\] 
where $\hu_{\text{potential}}$ is the parametrized potential and 
\[
\hu_{\text{BC}}(\bfx)=\prod_{i=1}^d\tanh\left(\frac{1}{2}x_i\right)\tanh\left(\frac{1}{2}(7-x_i)\right)\,,\qquad \bfx = [x_1, \dots, x_d]^T\,,
\] 
such that $\hu_{\text{BC}}(\bfx)=0$ and $\hu(\bfx; \bftheta(t)) = 0$ at the boundary of the domain $\Xcal$. The parametrization $\hu_{\text{potential}}$ of the potential is a two-layer fully connected network with sigmoid activation and 30 nodes per layer, with a total number of parameters $\nh = 1230$.  
As in previous examples, we use the RK4 scheme with a fixed time step size $\delta t=10^{-3}$.

\begin{figure}[t]
\centering
\begin{tabular}{cc}
    \resizebox{.48\linewidth}{!}{\large \input{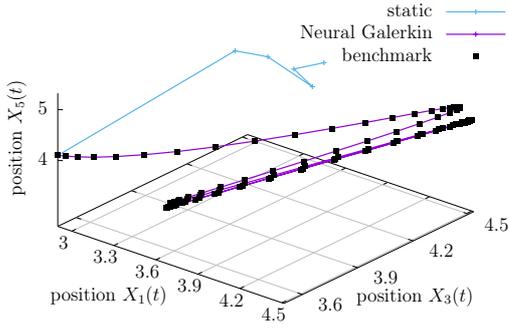}} & \resizebox{.48\linewidth}{!}{\large \input{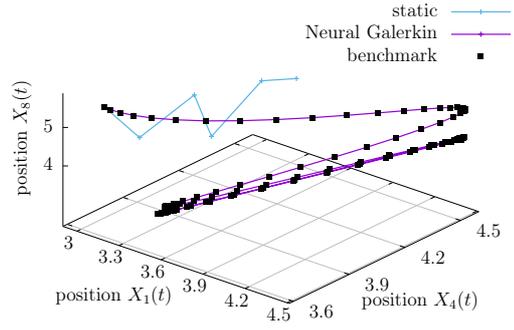}}\\
    \scriptsize (a) positions of physical particle $X_1(t), X_3(t), X_5(t)$& \scriptsize (b) positions of physical particle $X_1(t), X_4(t), X_8(t)$\\
    \resizebox{.48\linewidth}{!}{\large \input{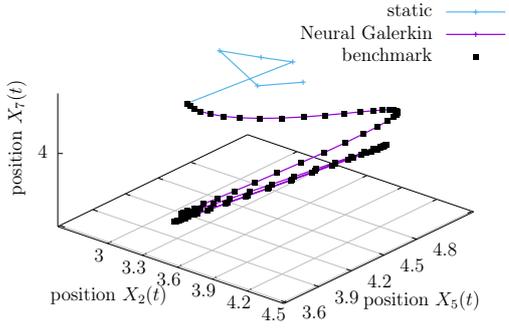}} & \resizebox{.48\linewidth}{!}{\large \input{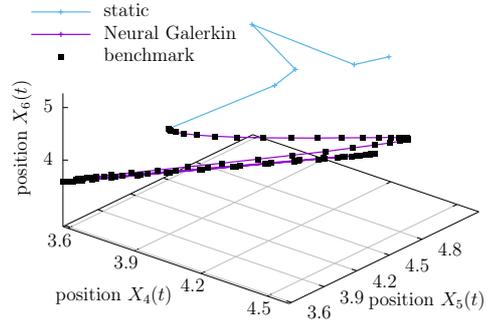}}\\
    \scriptsize (c) positions of physical particle $X_2(t), X_5(t), X_7(t)$ & \scriptsize (d) positions of physical particle $X_4(t), X_5(t), X_6(t)$
\end{tabular}
    \caption{The physical particles concentrate in the eight-dimensional domain over time, which means that the density function $u$ that is approximate by Neural Galerkin becomes local in the spatial domain. However, the dynamic particles allow Neural Galerkin schemes to adaptively keep track of the local behavior over time and give accurate approximations, whereas uniformly sampling the spatial domain leads to large errors in the predictions.}
    \label{fig:PTCubic8_3Dmean}
\end{figure}

\paragraph{Monte Carlo mean and covariance} Because no analytic solution is available but the mean and covariance can be estimated with Monte Carlo sampling, we evaluate the quality of numerical approximations by comparing their mean and covariance with those obtained with Monte Carlo sampling. The benchmark mean and covariance are estimated with Monte Carlo with $10^5$ paths from the SDE given in \eqref{eq:PT_particles}. The mean and covariance of the Neural Galerkin approximations $\hu$ are estimated with self-normalized importance sampling, where the biasing density is the Gaussian with the benchmark mean and covariance.

\begin{figure}[t]
\centering
\begin{tabular}{cc}
    \resizebox{.45\linewidth}{!}{\Large \input{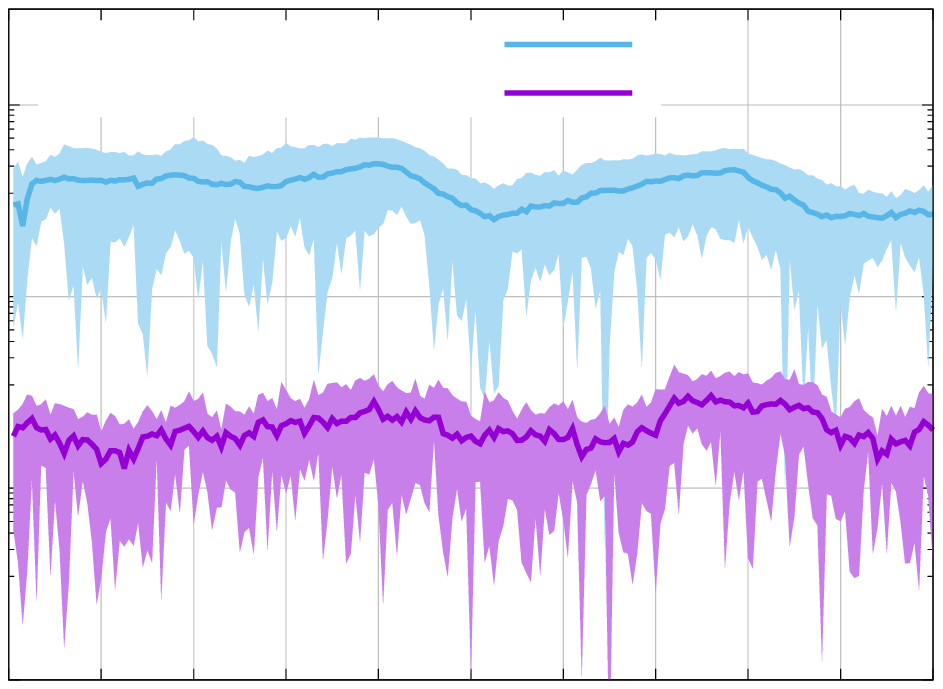}} & \resizebox{.45\linewidth}{!}{\Large \input{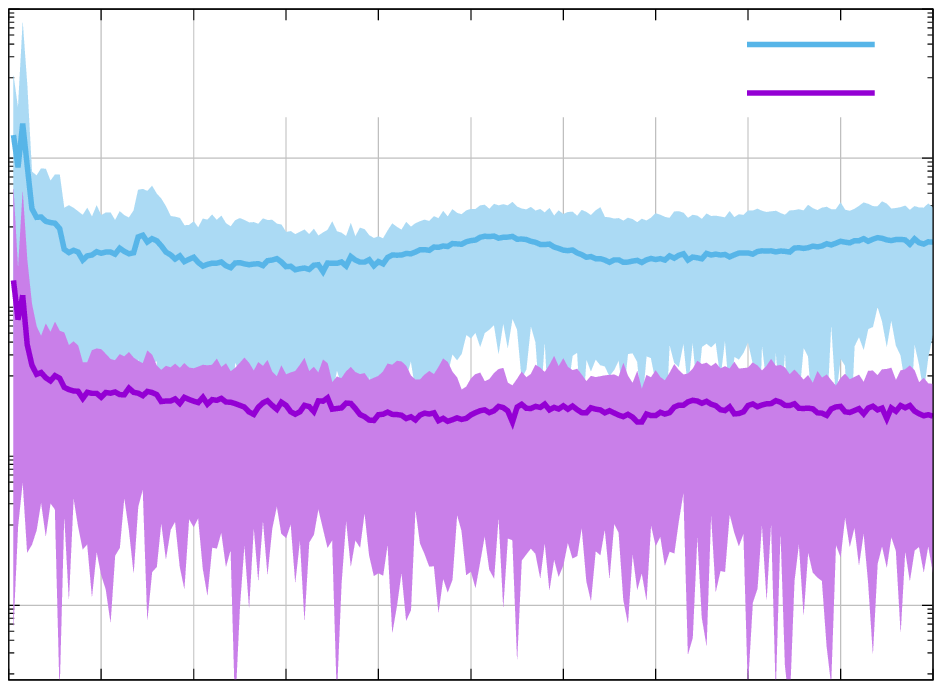}}\\
    \scriptsize (a) mean error & \scriptsize (b) covariance average error\\
    \resizebox{.45\linewidth}{!}{\Large \input{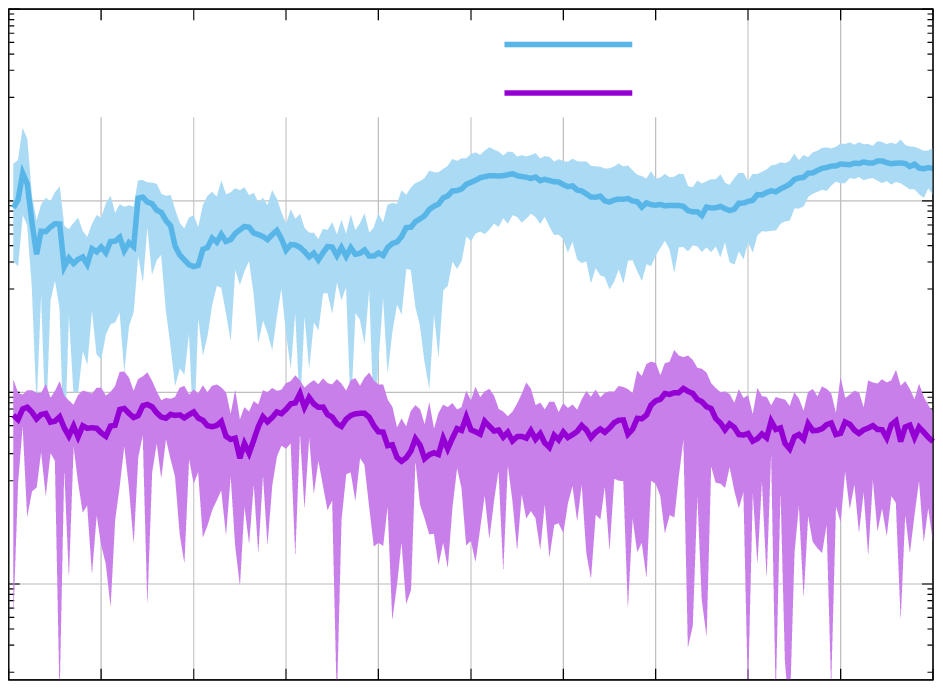}} & \resizebox{.45\linewidth}{!}{\Large \input{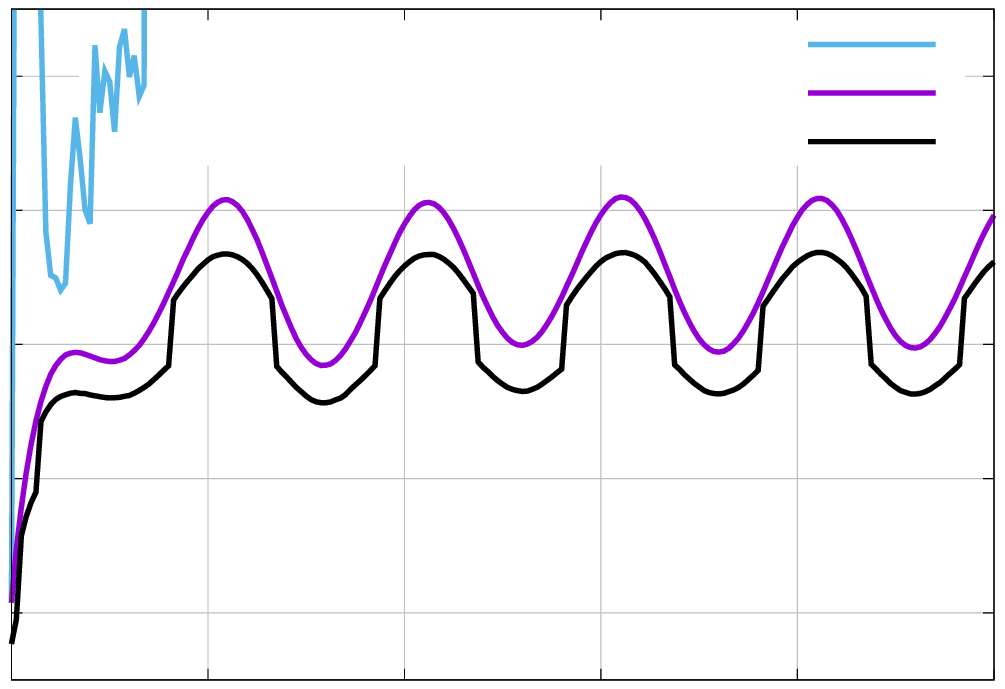}}\\
    \scriptsize (c) covariance error on diagonal & \scriptsize (d) entropy prediction
\end{tabular}
    \caption{With dynamic particles, Neural Galerkin schemes achieve orders of magnitude higher accuracy in predicting the mean  and covariance of the positions $X_1(t), \dots, X_d(t)$. The shaded region shows the minimum and maximum of the relative error over all dimensions, while the bold line shows the average relative error. Because Neural Galerkin schemes approximate the density, rather than providing sample paths as Monte Carlo approaches, quantities of interest that involve the density function such as the entropy can be computed.}
    \label{fig:PTCubic_err_plots}
\end{figure}

\paragraph{Quality of predicted mean and covariance} We compare the proposed Neural Galerkin scheme with dynamic particles to uniform sampling with $m = 2500$ samples over the eight-dimensional domain $\Xcal$. We first consider the case where the target measure $\mu_t$ is proportional to the squared residual. For the SVGD setup, we take 250 SVGD steps at each time step, with kernel bandwidth $0.05$ and SVGD step size $0.5$. The tempering parameter is set to $\gamma=0.5$.
As shown in Figure~\ref{fig:PTCubic8_3Dmean}, the mean of the positions of the physical particles described by \eqref{eq:PT_particles} concentrate over time so that the joint density $u$ becomes local in the eight-dimensional spatial domain. With $m=2500$ samples from the static, uniform measure, mean positions cannot be well approximated. In contrast, for the same number of $m = 2500$ samples, the proposed Neural Galerkin scheme provides accurate predictions of the positions. We quantify the accuracy in Figure~\ref{fig:PTCubic_err_plots} that shows the relative error in the mean and covariance with respect to the benchmark Monte Carlo results. Plot (a) of Figure~\ref{fig:PTCubic_err_plots} shows the relative error of the mean, which is orders of magnitude lower for Neural Galerkin approximations with dynamic particles than with uniform sampling. Similarly, the relative error of the covariance estimate obtained with Neural Galerkin and dynamic particles is also orders of magnitude lower than with uniform sampling. In fact, for uniform sampling, the relative error of the covariance is larger than one, which means the corresponding numerical solution is not predictive, which is in agreement with the results shown in Figure~\ref{fig:PTCubic8_3Dmean}. If we consider only the diagonal elements of the covariance matrix corresponding to the Neural Galerkin solution with dynamic particles, then the relative error is about one order of magnitude lower than the error averaged over all entries of the covariance matrix. 

\begin{figure}[t]
\centering
\begin{tabular}{cc}
    \resizebox{.45\linewidth}{!}{\Large \input{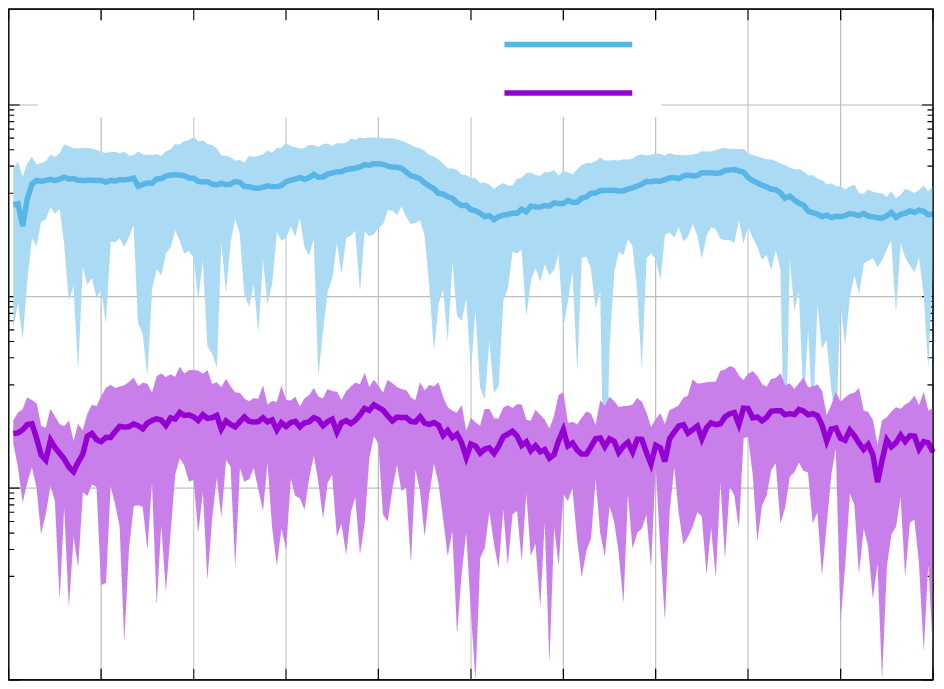}} & \resizebox{.45\linewidth}{!}{\Large \input{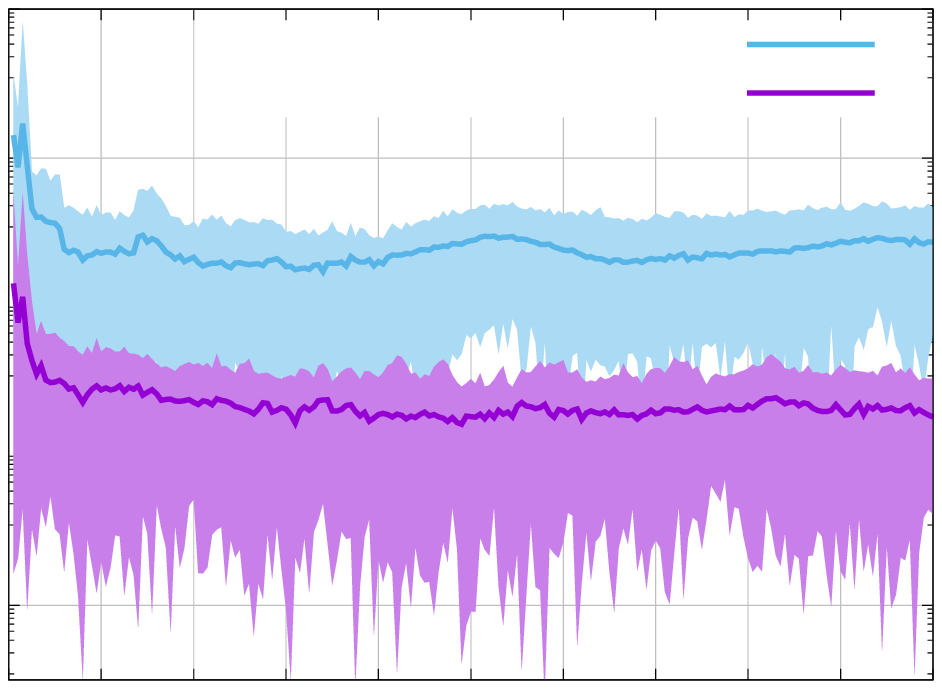}}\\
    \scriptsize (a) mean error & \scriptsize (b) covariance average error\\
    \resizebox{.45\linewidth}{!}{\Large \input{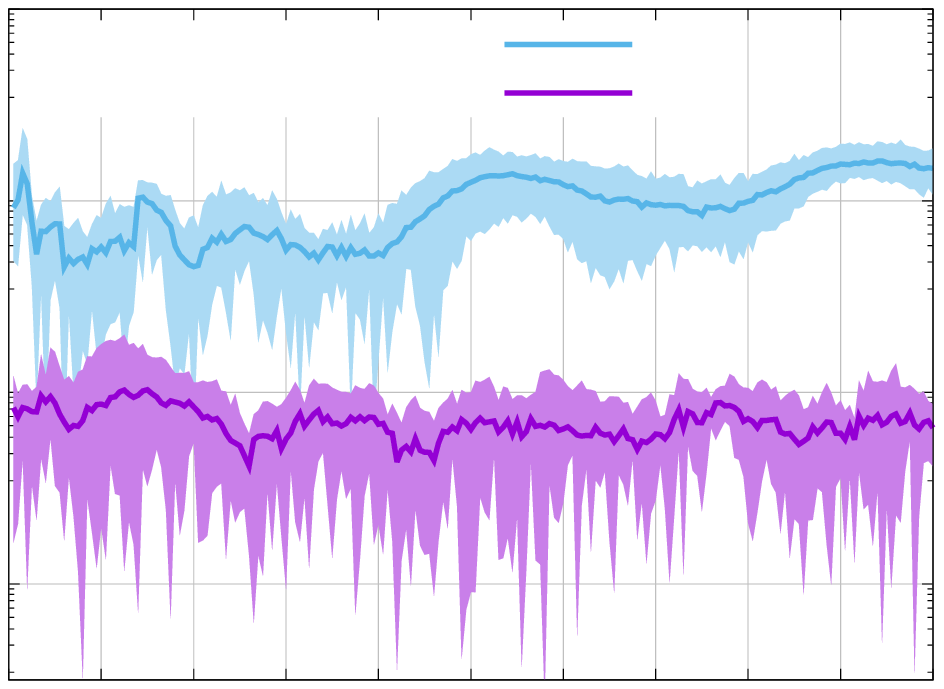}} & \resizebox{.45\linewidth}{!}{\Large \input{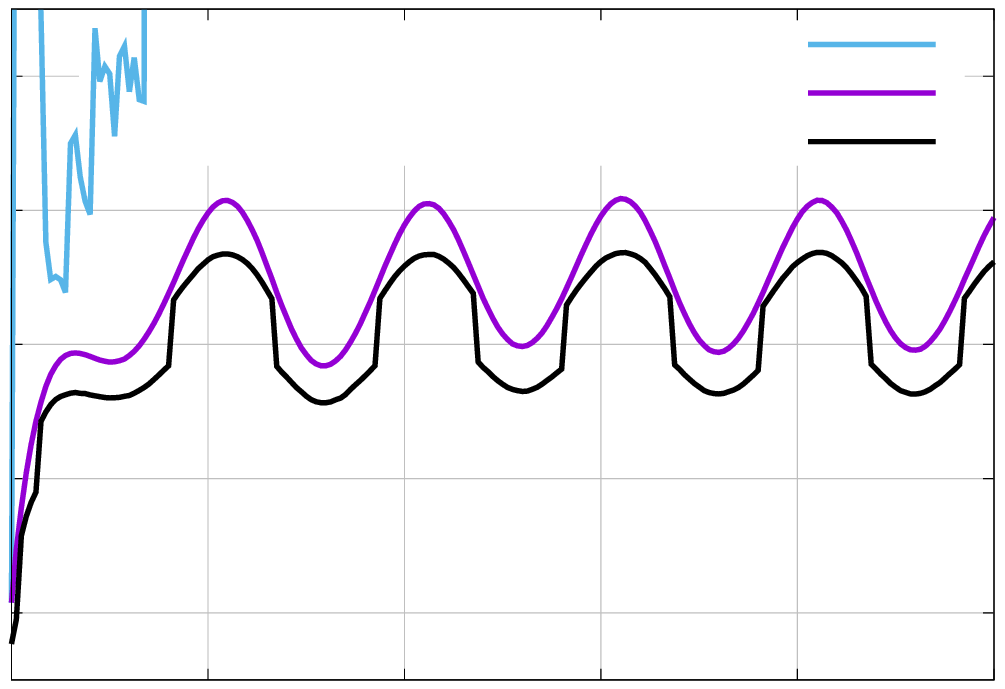}}\\
    \scriptsize (c) covariance error on diagonal & \scriptsize (d) entropy prediction
\end{tabular}
     \caption{This figure shows analogous results to Figure~\ref{fig:PTCubic_err_plots} except that in this figure the particle dynamics are induced by the magnitude of the solution as in \eqref{eq:NGDyn:PropMagnitudeSolution}.}
    \label{fig:PTCubic_SVGDsol_err_plots}
\end{figure}

\paragraph{Computing the entropy with Neural Galerkin schemes} In contrast to sampling paths of the SDE \eqref{eq:PT_particles} with Monte Carlo, we obtain an approximation of the density $u$ of the positions of $X_1(t), \dots, X_d(t)$ and thus can compute downstream quantities such as the entropy. For comparison purposes, we estimate a density from Monte Carlo paths of the SDE \eqref{eq:PT_particles} with kernel density estimation (KDE) and use it to compute an approximation of the entropy. Figure~\ref{fig:PTCubic_err_plots}(d) compares the entropy approximation obtained with Neural Galerkin schemes and the entropy obtained with KDE from Monte Carlo paths. For systems as the ones described by the SDE \eqref{eq:PT_particles}, it is known that the entropy oscillates smoothly, which is captured well by the Neural Galerkin approximation. The entropy approximation computed with Monte Carlo and KDE also captures the oscillations but is less smooth and shows spurious jumps.  

\paragraph{Sampling proportional to the magnitude of the numerical solution}
Now we impose dynamics on the particles with the target measure being proportional to the magnitude of the solution, as given in \eqref{eq:NGDyn:PropMagnitudeSolution}. In this case, we take 100 SVGD steps at each time step and set the kernel bandwidth to $5.0$ and the SVGD step size to $0.01$. The rest of the setup is the same as above: The tempering parameter is $\gamma=0.5$ and the number of particles is $m=2500$. Figure~\ref{fig:PTCubic_SVGDsol_err_plots} shows that the Neural Galerkin approximation with dynamic particles achieves a higher accuracy than the approximation obtained via uniform sampling.

\begin{figure}
\centering
\begin{tabular}{cc}
    \resizebox{.45\linewidth}{!}{\LARGE \input{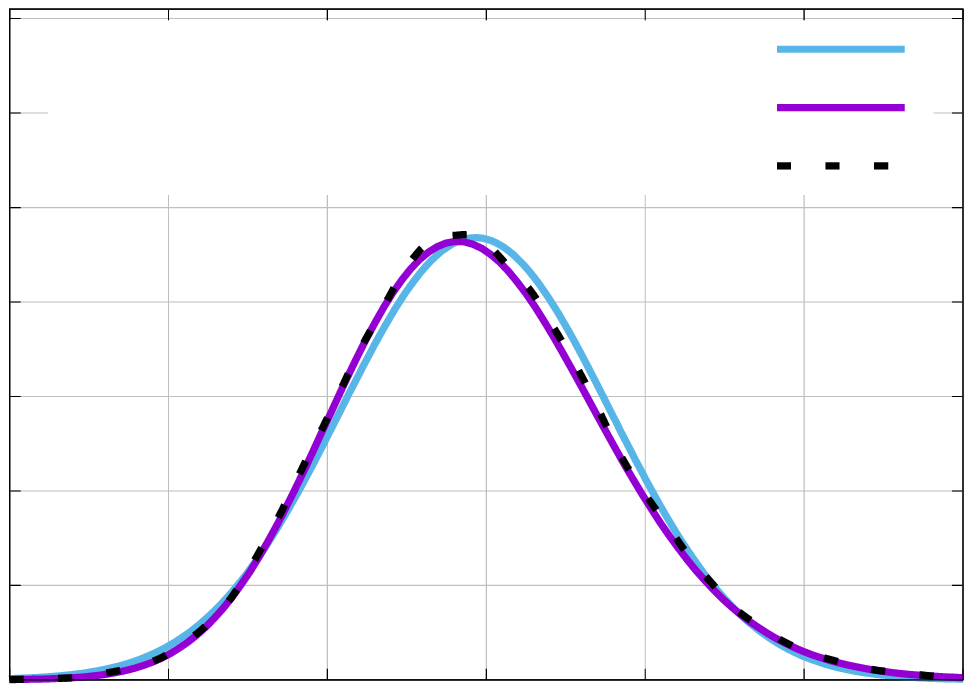}} & \resizebox{.45\linewidth}{!}{\LARGE \input{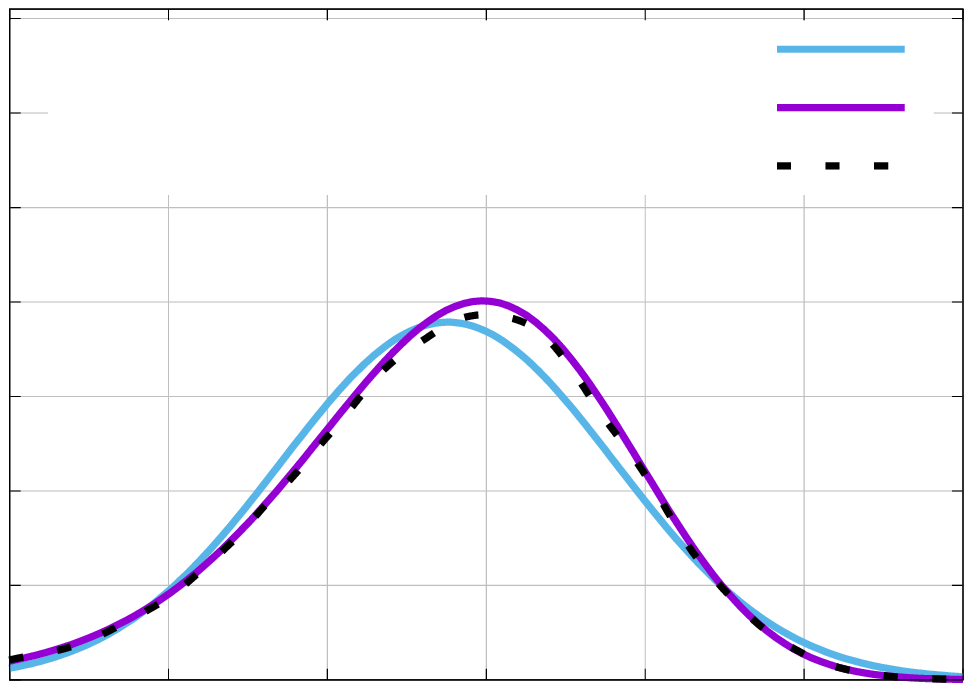}}\\
    \scriptsize (a) time $t=0.5$ & \scriptsize (b) time $t=1.5$\\
    \resizebox{.45\linewidth}{!}{\LARGE \input{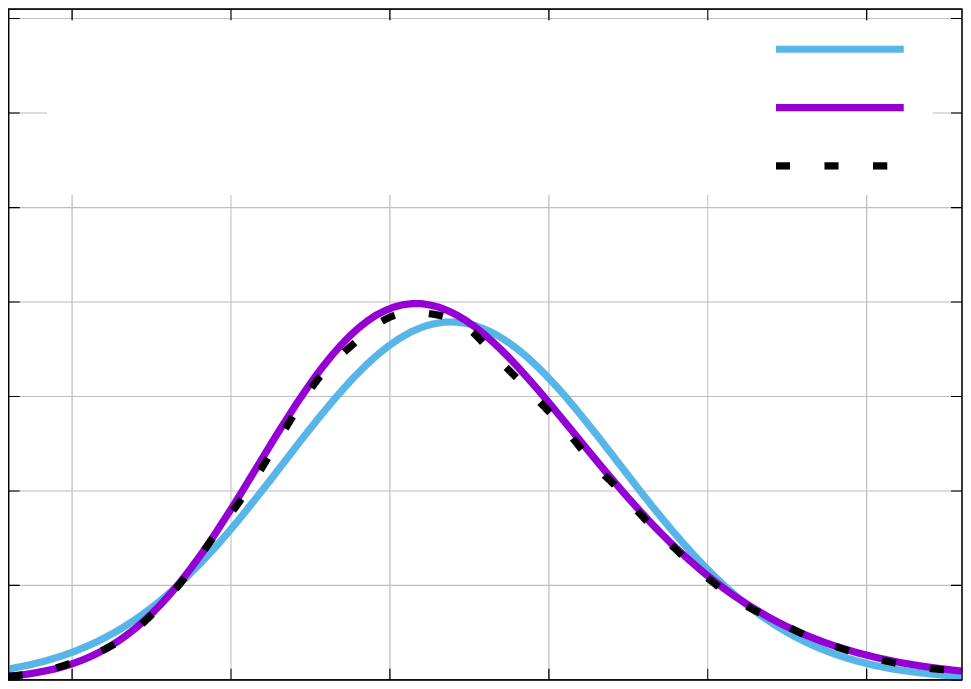}} & \resizebox{.45\linewidth}{!}{\LARGE \input{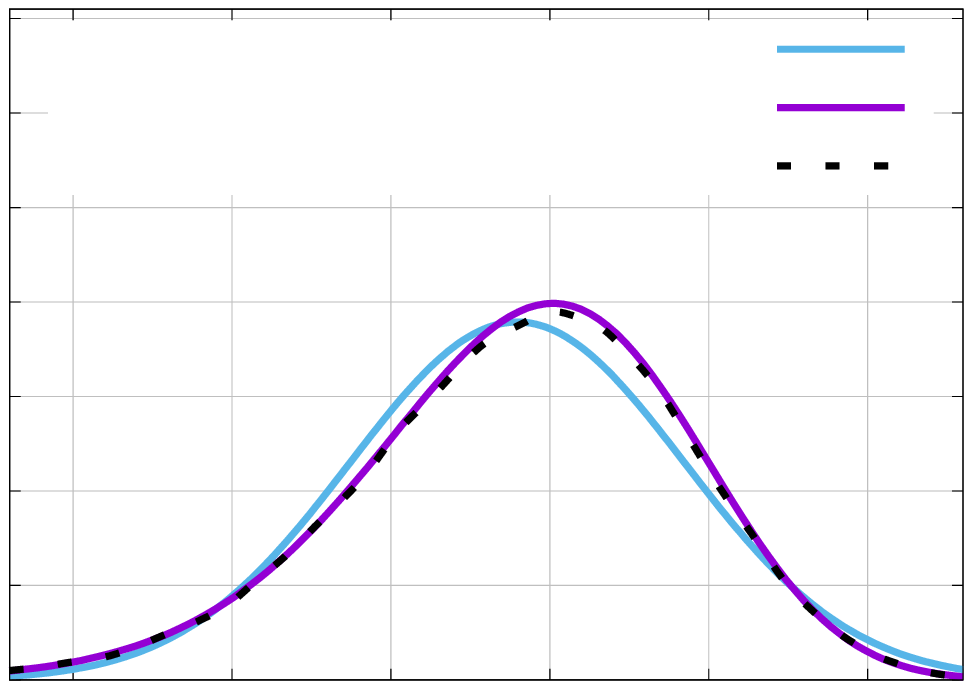}}\\
    \scriptsize (c) time $t=2.5$ & \scriptsize (d) time $t=3.5$\\
    \resizebox{.45\linewidth}{!}{\LARGE \input{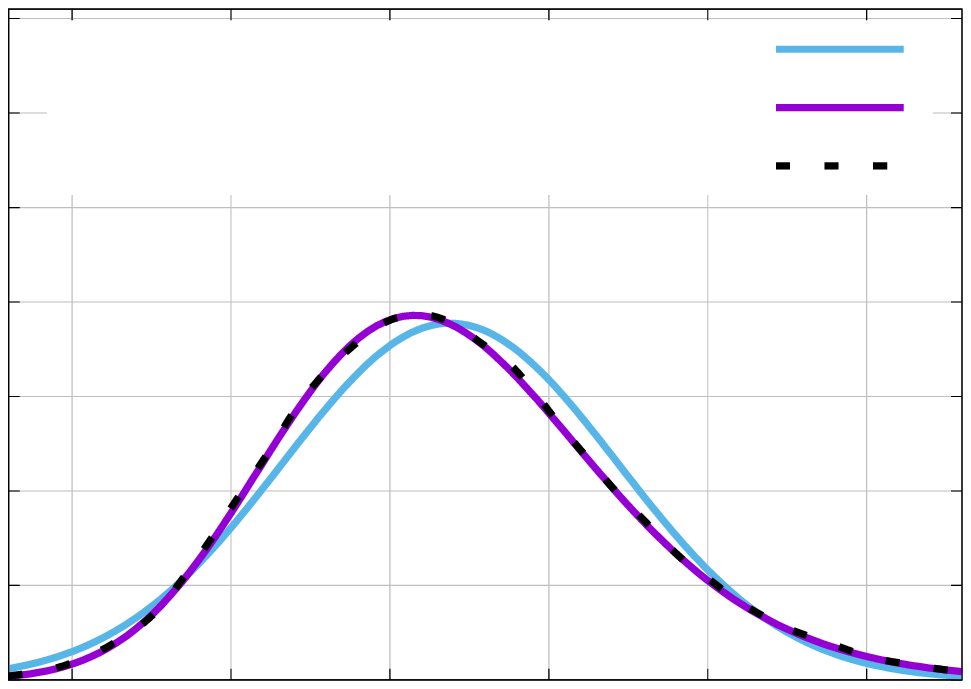}} & \resizebox{.45\linewidth}{!}{\LARGE \input{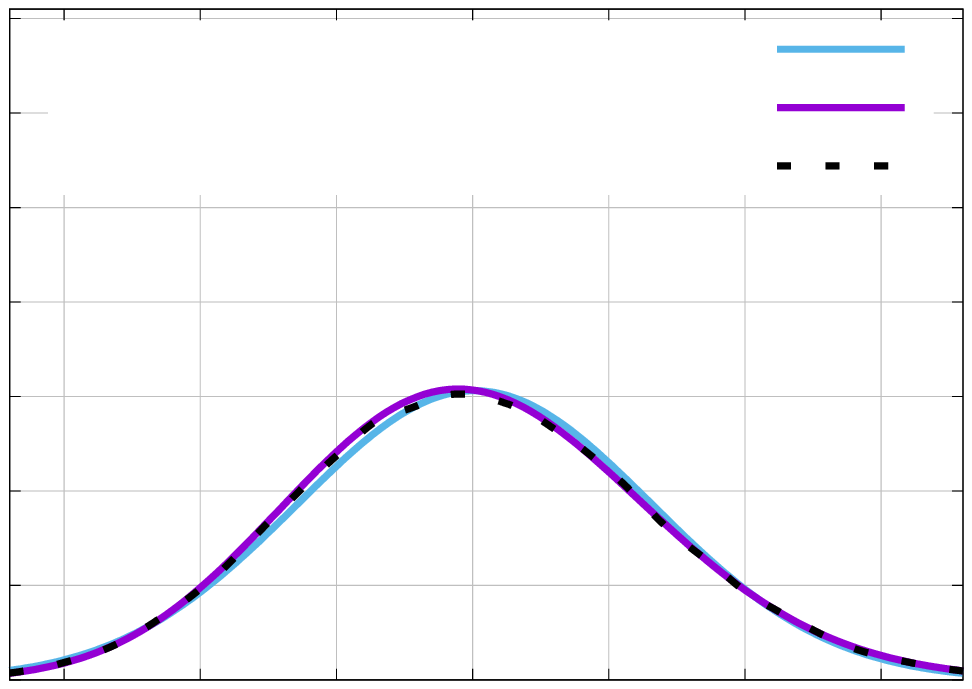}}\\
    \scriptsize (e) time $t=4.5$ & \scriptsize (f) time $t=5.0$\\
\end{tabular}
    \caption{The density $u$ given by the configuration of the Fokker-Planck equation that we consider is slightly different from a Gaussian density. This plot shows that Neural Galerkin approximations capture this difference well, which demonstrates the importance of approximating the density function rather than just the mean and covariance of the distribution.} 
    \label{fig:PTCubic8_SVGDsol_marginals}
\end{figure}

\paragraph{Capturing non-Gaussian behavior with Neural Galerkin schemes} Because the solution $u$ represents a density, one may be inclined to just approximate it with a Gaussian density, with mean and covariance that can be either estimated or computed via the solution of a system of ordinary differential equations; see \cite{NG22}. However, as we demonstrate in Figure~\ref{fig:PTCubic8_SVGDsol_marginals}, we consider a configuration of the Fokker-Planck equation that leads to non-Gaussian distributions that are captured well by the Neural Galerkin approximation. The results in Figure~\ref{fig:PTCubic8_SVGDsol_marginals} demonstrate the importance of approximating the density function rather than just the mean and covariance: The Neural Galerkin solution approximates well the density in general, rather than only the mean and covariance that would also be captured by a Gaussian approximation. 

\section{Conclusions}\label{sec:Conc}
Nonlinear parametrizations such as deep neural networks can achieve a faster error decay than traditional linear approximations from an approximation-theoretic perspective. However, numerically fitting the parameters to realize the fast error decay is challenging and critically depends on the available training data for estimating the population risk with the empirical risk. The results of this work indicate that actively adapting the measure that specifies where in the spatial domain to collect data for training is essential for numerically approximating well solution fields governed by evolution equations. The adaptive data collection is especially important when solution fields have local features that move through the spatial domain over time, such as wave fronts, phase transitions, and other coherent structures.

\section*{Acknowledgements}
The authors Wen and Peherstorfer were partially supported by the National Science Foundation under Grant No.~2046521 and the Office of Naval Research under award N00014-22-1-2728.

\bibliographystyle{abbrv}
\bibliography{main}

\end{document}